\documentclass[12pt,final,letterpaper]{article}
\usepackage{amsmath,amssymb,amsthm}
\newtheorem{thm}{Theorem}[section]
\newtheorem{prop}[thm]{Proposition}
\newtheorem{cor}[thm]{Corollary}
\newtheorem{lem}[thm]{Lemma}
\theoremstyle{definition}
\newtheorem{defn}{Definition}[section]
\newtheorem{rem}{Remark}[section]

\renewcommand{\H}{\mathbb{H}}
\newcommand{\N}{\mathbb{N}}
\newcommand{\R}{\mathbb{R}}

\newcommand{\al}{\alpha}
\newcommand{\be}{\beta}
\newcommand{\del}{\delta}
\newcommand{\eps}{\epsilon}
\newcommand{\ga}{\gamma}
\newcommand{\la}{\lambda}
\renewcommand{\phi}{\varphi}

\newcommand{\nul}{\emptyset}
\renewcommand{\O}{\Omega}
\newcommand{\sub}{\subseteq}
\renewcommand{\L}{\mathcal{L}^2}
\newcommand{\Haus}{\mathcal{H}^1}

\newcommand{\HB}{\mathrm{H}}
\def\Xint#1{\mathchoice
   {\XXint\displaystyle\textstyle{#1}}%
   {\XXint\textstyle\scriptstyle{#1}}%
   {\XXint\scriptstyle\scriptscriptstyle{#1}}%
   {\XXint\scriptscriptstyle\scriptscriptstyle{#1}}%
   \!\int}
\def\XXint#1#2#3{{\setbox0=\hbox{$#1{#2#3}{\int}$}
     \vcenter{\hbox{$#2#3$}}\kern-.5\wd0}}

\def\dashint{\Xint-}

\newcommand{\res}
{\mathop{\hbox{\vrule height 7pt width .5pt depth 0pt \vrule height .5pt width 6pt depth 0pt}}\nolimits}
\begin{document}
\newcommand{\Is}{\mathcal{I}}
\newcommand{\e}{\epsilon}
\newcommand{\p}{\partial}
\renewcommand{\d}{\delta}
\renewcommand{\a}{\alpha}
\renewcommand{\b}{\beta}
\newcommand{\s}{\sigma}
\newcommand{\z}{\zeta}

\newcommand{\B}{B_d}
\renewcommand{\P}{\mathcal P}
\newcommand{\V}{\mathcal V}
\renewcommand{\e}{\varepsilon}
\renewcommand{\eps}{\varepsilon}
\renewcommand{\theta}{\vartheta}
\renewcommand{\rho}{\varrho}
\renewcommand{\r}{\rho}
\newcommand{\LL}{\mathcal M}
\renewcommand{\S}{\mathcal S}
\newcommand{\wh}{\hat}
\newcommand{\GG}{\Gamma}
\title{Convex isoperimetric sets in the Heisenberg group}

\author{
   Roberto Monti\thanks{The first author would like to thank
                         the Institute of Mathematics of the
                         University of Bern, where part of
                         the work was done.} \\
   Dipartimento di Matematica Pura ed Applicata, \\
   Universit\`a di Padova, Via Belzoni 7, \\
   35100 Padova, Italy. \\
   email: monti@math.unipd.it
   \and
   Matthieu Rickly\thanks{The second author would like to thank
                           the Institute for Mathematical Research
                           of the ETH Z\"{u}rich, the Institute of
                           Mathematics of the University of Bern and
                           the Department for Pure and Applied Mathematics
                           of the University of Padova.} \\
   Institute of Mathematics, \\
   University of Bern, Sidlerstrasse 5, \\
   3012 Bern, Switzerland. \\
   email: matthieu.rickly@math.unibe.ch}



\date{July 26, 2006}

\maketitle

\begin{abstract}
We characterize convex isoperimetric sets in the Heisenberg group
endowed with horizontal perimeter. We first prove Sobolev
regularity for a certain class of $\R^2$-valued vector fields of
bounded variation in the plane related to the curvature equations.
Then, by an approximation-reparameterization argument, we show
that the boundary of convex isoperimetric sets is foliated by
geodesics of the Carnot-Carath\'eodory distance.
\end{abstract}


\maketitle

\section{Introduction}

We identify the Heisenberg group $\H$ with $\R^3 =\mathbb C
\times\R$ endowed with the group law
\begin{equation}
 \label{law!}
   (z,t)(z',t') = ( z+z', t+t'+ 2 \mathrm{Im} z\bar z'),
\end{equation}
where $t,t'\in\R$, $z=x+iy,z'=x'+iy'\in \mathbb C$ and
$\mathrm{Im}z\bar z' = y x'-x y'$. The Lie algebra of
left-invariant vector fields is spanned by
\begin{equation}
 \label{XY}
  X = \frac{\p}{\p x}+2y\frac{\p}{\p t}, \quad
  Y = \frac{\p}{\p y}-2x\frac{\p}{\p t} \quad
  \textrm{ and } \quad
  T = \frac{\p}{\p t},
\end{equation}
and the distribution spanned by $X$ and $Y$, called horizontal
distribution, generates  the Lie algebra by brackets.

The natural volume in $\H$ is given by the Haar measure, which, up
to a positive factor, coincides with Lebesgue measure $\mathcal L
^3$ in $\R^3$. Lebesgue measure can also be recovered as the
Riemannian volume induced by the left-invariant metric for which
$X$, $Y$ and $T$ are orthonormal. We denote the volume of a
measurable set $A \sub \H$ by $\V(A) =\mathcal L^3(A)$.

The horizontal perimeter (or simply perimeter) of a measurable set
$A\sub \H$ is
\begin{equation}
 \label{peri}
  \P(A) = \sup\left\{ \int _ A \big( X\phi_1+Y\phi_2 \big )\,
   d\mathcal L^3   \; \Big| \;
   \phi_1,\phi_2\in C_c^1(\R^3),\, \phi_1^2+\phi_2^2\leq
  1\right\}.
\end{equation}
If $\P(A)<+\infty$, the set $A$ is said to be of finite perimeter.
Perimeter is left-invariant and $3$-homogeneous with respect to
the group of dilations $\d_\la:\H\to\H$ given by
\begin{equation}
 \label{dila}
     \d_\la (z,t) = (\la z, \la^2 t), \quad \la>0,
\end{equation}
that is $\P(\d_\la(A))=\la^3 \P(A)$. Definition \eqref{peri} of
perimeter, which is modelled on De Giorgi's notion of perimeter
in Euclidean spaces, was introduced in \cite{GN} (see also \cite{FGW}).
If $A$ is smooth (e.g.~of class $C^2$), the perimeter coincides with the
Minkowski content and with the 3-dimensional Hausdorff measure,
precisely
\begin{equation}
  \label{titti}
   \P (A) = \lim_{\e\downarrow0}\frac{\V(A_\e)-\V(A)}{\e}
          = \mathcal S ^3(\p A).
\end{equation}
Here, $A_\e$ is the $\e$-neighbourhood of $A$ with respect to
Carnot-Carath\'eodory distance (we recall the definition in
Section \ref{prel}), and $\mathcal S^3$ is the spherical
$3$-Hausdorff measure, properly normalized, defined by means of
the same distance. The equality of perimeter and Minkowski
content is proved in \cite{MSC}, the equality of perimeter and
$\mathcal S^3$ is proved in \cite{FSSC}. Perimeter also admits an
integral-geometric representation formula (see \cite{Mo}).

Volume and perimeter are related via the isoperimetric
inequality
\begin{equation}
 \label{ISO}
         \V(A) \leq c \P(A)^{4/3},
\end{equation}
where $c>0$ is a constant and $A\sub\H$ is any measurable set with
finite perimeter and volume. This inequality is proved by P.~Pansu
in \cite{Pan1} and \cite{Pan2} for smooth domains, with $\mathcal
S ^3(\p A)$ instead of $\P(A)$. A set $A\sub\H$ with
$0<\V(A)<+\infty$ is called an isoperimetric set if it minimizes
the isoperimetric ratio
\begin{equation}
 \label{Is!}
    \Is(A) = \frac{\P(A)^{4/3}}{\V(A)}.
\end{equation}
The existence of isoperimetric sets is proved by Leonardi and
Rigot in \cite{LR} by a concentration-compactness
argument. Pansu notes that the boundary of a smooth isoperimetric
set has ``constant mean curvature'' and that a smooth surface has
``constant mean curvature'' if and only if it is foliated by
horizontal lifts of plane circles with constant radius. Then he
conjectures that an isoperimetric set is obtained by rotating around
the center of the group a geodesic joining two points in the center.
Recently, Pansu's conjecture reappeared in \cite{LM}.

Up to a left translation, the conjectured isoperimetric sets are
rotationally symmetric, i.e.~they are invariant with respect to
isometries which leave the center invariant. However, the non
commutative group law makes it difficult to prove by rearrangement
or symmetrization that the isoperimetric ratio \eqref{Is!} is minimized
by rotationally symmetric sets.

Assuming rotational symmetry, it is easy to determine the isoperimetric
profile, as shown in \cite{M} and \cite{DGN} in some special
cases, and in \cite{RR1} in the general case. In fact, it suffices
to assume the rotational symmetry of a certain horizontal section (see
\cite{GN2} and \cite{R}). Isoperimetric sets can be also determined
assuming regularity instead of symmetry. This is an important new result
due to M.~Ritor\'e and C.~Rosales. In \cite{RR2}, they prove Pansu's conjecture
for isoperimetric sets of class $C^2$, without any symmetry assumption.
However, a regularity theory for isoperimetric sets in the Heisenberg group
does not yet exist (but see \cite{CHMY} and \cite{Pau2}).

It is worth mentioning that a 2-dimensional version of the problem is
formulated and solved by the first author and Morbidelli in
\cite{MM}. In this setting, isoperimetric sets coincide with the
section of the set \eqref{bubble} below with the $y=0$ plane,
properly scaled and translated. Finally, we refer the reader to the
monograph \cite{CDST} for a more detailed introduction to the
isoperimetric problem in the Heisenberg group.

\bigskip

In this article, we characterize isoperimetric sets in the class of convex
sets. A set $C \sub \H$ is convex if it is convex with respect to the
usual vector space structure of $\R^3$. We do not assume any symmetry
or regularity besides convexity.

\begin{thm}[Convex isoperimetric sets]
  \label{UN}
Up to a left translation and a dilation, any convex open
isoperimetric set in $\H$ coincides with the set
\begin{equation}
 \label{bubble}
   \left\{ (z,t)\in\H \;\big|\; |t| <
   \arccos|z|+|z|\sqrt{1-|z|^2}, \, |z|<1\right\}.
\end{equation}
\end{thm}

\noindent The set in \eqref{bubble} is foliated by Heisenberg geodesics
(see Section \ref{prel}), it is globally of class $C^2$, but it fails to
be of class $C^3$ at the north and south poles $(0,\pm \pi/2)$. At these
points, the plane spanned by the vector fields $X$ and $Y$, the horizontal
plane, is tangent to the boundary of the set.

In the first part of the paper, we describe the characteristic set
of a convex set $C\sub \H$. We say that a point $p\in \p C$ is
characteristic if the horizontal plane at $p$ is a supporting
plane for $C$ at $p$ and we write $p \in \Sigma(C)$. We show that
$\Sigma(C) = \Sigma^- \cup \Sigma^+$, where $\Sigma^-, \, \Sigma^+$
are two closed, disjoint, horizontal segments (possibly points).
This is Theorem \ref{H6} in Section \ref{konvex}.

In Section 4, we derive the curvature equations by a variational argument.
We consider a convex set
\begin{equation}
 \label{dado!}
  C = \big\{ (z,t) \in \H \;|\; z \in D, \, f(z) \leq t \leq g(z) \big\},
\end{equation}
where $D\sub\R^2$ is a closed convex set in the plane with nonempty
interior, and $-g,f:D \to \R$ are convex functions. The intersection of
$\Sigma(C)$ with the graph of $f$ is a union of two (possibly empty) line
segments. We denote its projection onto the $xy$-plane by $\Sigma(f)$.

If $C \sub \H$ is isoperimetric, then the function $f:D \to \R$
satisfies the partial differential equation
\begin{equation}
   \label{distr!}
   \mathrm{div} \left( \frac{\nabla f(z)+2z^\perp}
   {|\nabla f(z)+2z^\perp|} \right) = \frac{3\P(C)}{4\V(C)}
\end{equation}
in $\mathrm{int}(D) \setminus \Sigma(f)$. Here $z = (x,y)$ and
$z^\perp = (-y,x)$. The curvature operator in \eqref{distr!} has
been studied by several authors (besides the previous references,
see also e.g.~\cite{Pau1}, \cite{Pau2}, \cite{CHY}). The number
\begin{equation}
 \label{H!}
     H = \frac{3\P(C)}{4\V(C)}
\end{equation}
is the (horizontal) curvature of $\p C$, and the equation states
that the boundary of a convex isoperimetric set has constant
curvature in a weak sense.

If we write $u (z)= \nabla f(z)+2z^\perp$, then we can interpret
equation \eqref{distr!} as an equation for the distributional
derivative of $u/|u|$, which is a measure. In Section \ref{regular!},
we prove that this measure is absolutely continuous. In other words,
we get an improved regularity along horizontal directions for the
boundary of convex isoperimetric sets. This result is a corollary
of the following regularity theorem for $BV$ vector fields.

\begin{thm}[Improved regularity]\label{QUA}
Let $\O \sub \R^2$ be a bounded open set and let $u \in
BV(\O;\R^2)$ be a vector field. Suppose that:
\begin{enumerate}
\item
There exists $\del > 0$ such that $|u(z)| \geq \del$ for
$\L$-a.e.~$z \in \O$;
\item
$\mathrm{div}\, u^\perp \in L^1(\O)$;
\item
$\mathrm{div} \Big(\displaystyle\frac{u}{|u|}\Big) \in L^1(\O)$.
\end{enumerate}
Then $u/|u| \in W^{1,1}(\O;\R^2)$.
\end{thm}

In the reconstruction argument for isoperimetric sets in Section
\ref{final}, we also need to work with graphs of the form $x =
h(y,t)$, where $h:E \to\R$ is a convex function defined on some
bounded convex set $E \sub \R^2$ with nonempty interior. For these
graphs, the partial differential operator related to the
horizontal curvature is a bit more complicated than the one
appearing in \eqref{distr!} (see equation \eqref{pox}). In Theorem
\ref{Reggae}, we prove a more general version of Theorem
\ref{QUA}, which includes both curvature operators above as
special cases.

The improved regularity of the boundary is the starting point for
the geometric characterization of convex isoperimetric sets. A
first interesting consequence of Theorem \ref{QUA} is the
stability of the curvature with respect to smooth approximations,
which is proved in Section \ref{letto}.

\begin{thm}
  [Stability of the curvature]
  \label{CINQ}
Let $\O \sub \R^2$ be a bounded open set, let $f \in \mathrm{Lip}(\O)$
be a Lipschitz function such that $\nabla f \in BV(\O;\R^2)$ and denote by
$f_{\e} \in C^\infty(\O_\e)$  the standard mollification of $f$, with
$\O_\e=\{z\in\O\;|\; \mathrm{dist}(z,\p\O)>\e\}$. We consider the vector
fields
\begin{equation}
 \label{ueu!}
  u(z) = \nabla f(z) + 2z^\perp
    \quad\mbox{and}\quad
  u_\e(z) = \nabla f_\e(z) + 2z^\perp.
\end{equation}
Assume that:
\begin{enumerate}
\item
  There exist $\del>0$ and $\e_0>0$ such that
  $|u(z)|\geq \del$ for $\L$-a.e.~$z\in \O$
  and $|u_\eps(z)|\geq \del$ for all $z\in \O_\eps$,
  $0< \e <\e_0$;
\item
   $\displaystyle \mathrm{div} \Big(\frac{u}{|u|}\Big)\in L^1(\O)$.
\end{enumerate}
Then
\begin{equation}
   \label{ELLE!}
      \lim_{\eps\to0} \int_K \Big| \mathrm{div}
      \Big(\frac{u_\eps}{|u_\eps|}\Big)-
      \mathrm{div}\Big(\frac{u}{|u|}\Big)\Big|\,
      d z  = 0
\end{equation}
for any compact set $K \sub \O$.
\end{thm}

The vector field $v(z) = - u^\perp(z) =  2z-\nabla f^\perp(z)$
belongs to $BV_{\mathrm{loc}}(\mathrm{int}(D);\R^2)$ if $f:D\to\R$ is
convex. Moreover, its distributional divergence is in $L^\infty$,
in fact
\begin{equation}
 \label{div!}
   \mathrm{div}\,v  = 4
   \quad \mbox {in $\mathrm{int}(D)$}.
\end{equation}
Therefore, we can use the theory on the Cauchy Problem for vector
fields of bounded variation recently developed by Ambrosio in
\cite{A}. This theory extends the work of DiPerna and Lions on the
flow of Sobolev vector fields (\cite{DL}) to the $BV$ setting. The
bound on the divergence ensures the existence of a unique regular
Lagrangian flow $\Phi:K\times[-\r,\r] \to D$ ($\r > 0$ small
enough) relative to $v$ starting from a compact set $K \sub
\mathrm{int}(D)$. For $\L$-a.e.~$z \in K$, the curve $s \mapsto
\Phi(z,s)$ is an integral curve of $v$ passing through $z$ at time
$s = 0$.

In Sections \ref{Lag}, \ref{xillo} and \ref{Sobolev}, we show
that $\L$-a.e.~integral curve of $v$ is an arc of circle.
Since the horizontal lift of the flow $\Phi$ foliates the graph of the
function $f$, it follows that the boundary of the convex isoperimetric
set is foliated by geodesics.

\begin{thm}[Foliation by circles]
 \label{SIX}
Let $C\sub\H$ be a convex isoperimetric set with curvature $H>0$.
Then, for $\L$-a.e.~$z\in K$, the integral curve $s\mapsto
\Phi(z,s)$ is an arc of circle with radius $1/H$.
\end{thm}

The precise version of this result is Theorem \ref{effe} in
Section \ref{Lag}. In Theorem \ref{acca} we prove an analogous
result for graphs of the form $x = h(y,t)$. We give two different
proofs of these theorems. The proof in Section \ref{xillo} is
based on the stability property of regular Lagrangian flows and on
the stability property \eqref{ELLE!} of the curvature. Here, we
describe the proof contained in Section \ref{Sobolev}, which is
based on the Sobolev regularity for $u/|u|$ of Theorem \ref{QUA}
and on the regularity of the Lagrangian flow.

The vector field $v$ has a regular Lagrangian flow $\Phi:K\times[-\r,\r]\to
D$ but is only in $BV_{\mathrm{loc}}(\mathrm{int}(D);\R^2)$. On the
other hand, $v/|v|$ is in $W^{1,1}_{\mathrm{loc}}(\mathrm{int}(D);\R^2)$,
but its divergence is only in $ L^1_{\mathrm{loc}}(\mathrm{int}(D))$. In
order to compute the second order derivative of a generic integral curve
of $v$, we introduce suitable reparameterizations $\ga (s)= \Phi(z,\tau(s))$.
In Theorem \ref{regularity}, we show that for any vector field $w\in
W^{1,1}(\O;\R^2)$ defined in some open neighbourhood $\O$ of $K$, the curve
$\kappa(s)= w(\ga(s))$ is in $W^{1,1}$ (for $\L$-a.e.~$z\in K$) and
moreover
\begin{equation}
  \label{chain!}
   \dot\kappa = (\nabla w \circ \gamma) \, \dot{\ga}
   \quad \mbox{in the weak sense.}
\end{equation}

\noindent Using the chain rule \eqref{chain!} with $w = v/|v|$,
we can give a pointwise meaning to equation \eqref{distr!}
along the flow, and the integral curves of $v$ turn out to be arcs
of circles. The proof of Theorem \ref{acca} given in Section
\ref{Sobolev} is based on the same technique.

We combine Theorem \ref{effe} and Theorem \ref{acca} in Section
\ref{final}, where we finish the proof of Theorem \ref{UN}.
\section{Preliminaries}

\label{prel}
\setcounter{equation}{0}

We can identify the horizontal plane spanned by the vector fields
$X$ and $Y$ at the point $p=0$ with the $xy$-plane
\[
   \mathrm{H}_0 = \big\{ (z,0)\in\H \; | \; z\in\mathbb C \big\}.
\]
The $yt$-plane is the subgroup $ \{ (z,t)\in \H \; | \;
\mathrm{Re}\, z=0  \}$. In general, if $p = (z,t) \in \H$, we
define the horizontal plane at $p$ as the left translation of
$\mathrm H_0$ with $p$, that is
\[
   \mathrm{H}_p
   =
   p \mathrm{H}_0
   =
   \big\{ (z+z',t + 2\mathrm {Im} z\bar z' ) \in \H  \; | \; z' \in \mathbb C
   \big\}.
\]
The base point $p$ is uniquely determined by $\mathrm{H}_p$. The plane
$\mathrm{H}_p$ is the boundary of the two halfspaces
\begin{equation}
 \label{HS}
 \begin{split}
   \mathrm{H}_p^-  & = \big\{ (z',t')\in\H  \; | \;
      t' \leq t + 2 \mathrm{Im} z\bar z' \big \},
          \\
   \mathrm{H}_p^+ & = \big\{ (z',t')\in\H  \; | \;
      t' \geq  t + 2 \mathrm{Im} z\bar z' \big \}.
\end{split}
\end{equation}

We equip the tangent bundle $\mathrm{T}\H$ with the left invariant
inner product $\langle \cdot \, , \cdot \rangle$ which makes $X$,
$Y$ and $T$ orthonormal. An absolutely continuous path $\ga:[0,1] \to \H$
is said to be horizontal if $\dot{\ga}(s) \in \mathrm{span}_{\R}
\{ X(\ga(s)),Y(\ga(s)) \}$ for a.e.~$s \in [0,1]$. The curve $\ga
= (\ga_1,\ga_2,\ga_3)$ is horizontal if and only if
\begin{equation}
 \label{lift}
 \ga_3 (s) = \ga_3(0) + 2 \int_0^s \big ( \dot\ga_1\ga_2 - \dot\ga_2
 \ga_1 \big)\, d\s.
\end{equation}
We call the plane curve $\kappa = (\ga_1,\ga_2)$ the horizontal
projection of $\ga$. If $\kappa=(\ga_1,\ga_2)$ is a given plane
curve, then a curve $\ga=(\kappa,\ga_3)$  with $\ga_3$ given by
\eqref{lift} for some $\ga_3(0)$ is called a horizontal lift of
$\kappa$.

The sub-Riemannian length of $\ga$ is
\begin{equation}
 \label{lella}
   \mathrm{L}(\ga) = \int_{0}^{1} \langle
   \dot{\ga}(s) , \dot{\ga}(s) \rangle^{1/2} \, ds
   = \int_0^1|\dot\kappa(s)|\, ds,
\end{equation}
where $\kappa$ is the horizontal projection of $\gamma$ and
$|\dot\kappa|$ is the Euclidean length of $\dot\kappa \in\R^2$. A
path $\ga:[0,1] \to \H$ is said to be admissible for a given pair
$(p_0,p_1) \in \H \times \H$ if $\ga(0) = p_0$, $\ga(1) = p_1$,
and $\ga$ is an absolutely continuous, horizontal curve. The
Carnot--Carath\'{e}odory distance of points $p_0, p_1 \in \H$ is
\[
   d(p_0,p_1)
   =
   \inf \big\{ \mathrm{L}(\ga) \; | \; \mbox{$\ga$
      is admissible for $(p_0,p_1)$} \big\}.
\]
The open Carnot--Carath\'{e}odory ball of radius $r$ centered at
$p \in \H$ is denoted by $B_d(p,r)$. The distance $d$ induces the
Euclidean topology on $\H$. The metric space $(\H, d)$ is
complete, locally compact and geodesic.

Geodesic curves (i.e.~lenght minimizing curves between points) can be
computed explicitly, by minimizing the length functional on the right
hand side of \eqref{lella} subject to the constraint \eqref{lift} with
$s = 1$.  Take points $p_0 = 0$ and $p_1 = (z,t) \in \H$. We have the
following cases:
\begin{itemize}
  \item[1)] If $t=0$, the geodesic $\ga:[0,1]\to\H$ connecting $0$ to $p_1=(z,0)$
  is the line segment $\ga(s) = (s z,0)$, $s\in[0,1]$.

  \item[2)] If $t>0$ and $z\neq 0$, the geodesic connecting $0$
  to $p_1=(z,t)$ is the horizontal lift (starting at $0$)
  of the arc of circle from $0$ to $z$ in the $xy$-plane
  (oriented clockwise), such that the plane region  bounded by
  the arc and by the segment joining $0$ to $z$ has area equal to
  $t/4$. This geodesic is unique.

  \item[3)] If $t>0$ and $z=0$, the geodesic from $0$ to
  $p_1=(0,t)$ is not unique. Take any full circle (oriented
  clockwise) passing through $0$ and with area equal to $t/4$.
  The horizontal lift of the circle (starting from $0$) is the
  desired geodesic.
\end{itemize}

\noindent The case $t<0$ is similar and the case $p_0\neq0$ is
obtained by left translation. If the arc of circle in 2) and 3)
has radius $0 < R < +\infty$, we say that the geodesic has
curvature $H = 1/R$.

The union of all geodesics joining $(0,-\pi/2) $ to $(0,\pi/2)$ of
case 3) bounds the isoperimetric set conjectured by Pansu. The
horizontal lift of the plane circle $\kappa(s) = \frac 12 \big(
1+\cos s ,- \sin s\big)$,
  $s\in[-\pi,\pi]$,
passing through the point $(1,0,0)\in\H$ at time $s=0$ is the
curve $\ga:[-\pi,\pi]\to\H$
\begin{equation}
  \label{seed}
     \ga(s) = \frac 12 \big( 1+\cos s, -\sin s, s+\sin s \big).
\end{equation}
The third coordinate can be computed using formula \eqref{lift}.
The curve $\ga$ is a geodesic with curvature $H=2$, starting from
$\ga(-\pi) = (0,0,-\pi/2)$ and reaching $\ga(\pi) = (0,0,\pi/2)$. If
$(z,t) =\ga(s) \in \ga([-\pi,\pi])$ is a point on the curve, then
we have
\[
   |z| = \Big( \frac{1+\cos s}{2} \Big)^{1/2}\quad\textrm{ and
   }\quad
   t = \frac{1}{2}\big(s+\sin s\big),
\]
and we obtain the relation $|t| = \arccos|z|+|z|\sqrt{1-|z|^2}$.
We call the set
\begin{equation}
 \label{bubble2}
   C_{2} = \left\{ (z,t)\in\H \;\big|\; |t| <
   \arccos|z|+|z|\sqrt{1-|z|^2}, \, |z|<1\right\}
\end{equation}
isoperimetric bubble with curvature $H=2$. The boundary of $C_{2}$
is a compact surface which is obtained by rotation of the generating
curve \eqref{seed} around the $t$-axis. This surface is globally of
class $C^2$, but fails to be of class $C^3$ at the characteristic
points $(0,\pm \pi/2)$.

We conclude this preliminary section with the following
representation formula for perimeter. If $A \sub \H$ is a bounded
open set such that  $\partial A$ is a Lipschitz surface, then
\begin{equation}
  \label{area}
   \P(A) = \int_{\partial A} \sqrt {( X \cdot \nu)^2 + (Y \cdot \nu)^2} \, d \mathcal H^2,
\end{equation}
where $\nu$ is a unit normal to $\partial A$. Here and in the following,
$\cdot$ denotes the standard inner product in $\R^3$ or $\R^2$
(we think of $X$ and $Y$ as vectors in $\R^3$). $\mathcal H^2$ is
the 2-dimensional Hausdorff measure in $\R^3$ (with respect to the
usual Euclidean distance). For a proof of formula \eqref{area}, see
\cite{FSSC}. 
\section{Characteristic set of convex sets and convex functions}

\label{konvex}

\setcounter{equation}{0}

A set $C \sub \H$ is convex if it is convex with respect to the
standard convexity of $\R^3$ as a vector space. Note that this
notion is invariant with respect to affine transformations of
$\R^3$. In particular, it is invariant with respect to left
translations, dilations and isometries of $\left( \H, d \right)$.
Given $p \in \partial C$ and a plane $\pi \sub \H$ passing through $p$, we
say that $\pi$ is a supporting plane for $C$ at $p$ if $\pi \cap
\overline{C} \sub \partial C$. The property of being a supporting
plane for some convex set $C$ is a local property.  The
characteristic set of $C$ is
\[
   \Sigma (C)
   =
   \big\{
      p \in \partial C \; | \;
      \mbox{$\mathrm{H}_p$ is a supporting plane for $C$ at $p$}
   \big\},
\]
where $\mathrm{H}_p$ is the horizontal plane with base point $p$.

\begin{lem}
 \label{approx}
Let $C \sub  \H$ be a convex set such that $\B(p,r_1) \sub C \sub
\B(p,r_2)$ for some $p \in \H$ and suitable $0 < r_1 < r_2$. Then
there exists a sequence $\{ C_k \}_{k \in \N}$ of strictly convex
sets of class $C^\infty$ such that $\B(p,r_1/2) \sub C_k \sub
\B(p,2r_2)$ for all sufficiently large $k \in \N$ and $C_k \to C$
as $k\to\infty$ with respect to the Hausdorff distance.
\end{lem}

\begin{proof}
After a left translation, we can assume that $p = 0$. Consider the
functions $f_k = \chi_{1/k} \ast (g + | \cdot |^2/k)$, where
$\chi$ is a smoothing kernel and $g$ is the homogeneous convex
gauge function associated to $C$. Then each $f_k$ is a smooth,
strictly convex function, and the sets $C_k = \left\{ p \in \R^3
\; | \; f_k(p) \leq 1 \right\}$ have the required properties.
\end{proof}

\begin{lem}\label{lello}
If $C \sub \H$ is a convex, bounded open set with $C^\infty$
boundary, then $\Sigma (C) \neq \nul$.
\end{lem}

\begin{proof}
For $p \in \partial C$, let us denote by $V(p) \in \mathrm{T}_p
\partial C$ the orthogonal projection of the vector field $T$
onto $\mathrm{T}_p \partial C$ with respect to the left invariant
inner product $\langle \cdot \, , \cdot \rangle$. Since $\partial C$
is diffeomorphic with the sphere $S^2$, there exists
$p \in \partial C$ such that $V(p) = 0$. Hence $\mathrm{T}_p
\partial C = \mathrm{H}_p$, and therefore $p \in \Sigma (C)$.
\end{proof}

\begin{lem}
  \label{lallo}
If $C \sub \H$ is a strictly convex, bounded set, then
$\Sigma (C)$ contains at most two points.
\end{lem}

\begin{proof}
Let $p_0, \, p_1 \in \Sigma (C)$ be two characteristic points such
that  $p_0 \neq p_1$. After a left translation, we can assume that
$p_0 = 0$ and $p_1= (z_1,t_1)$ with $t_1 > 0$. The convex set $C$
is contained between the horizontal planes $\mathrm{H}_{p_0}$ and
$\mathrm{H}_{p_1}$. By strict convexity of $C$, the relative
interior of the line segment $S = \left\{ (s z_1, s t_1) \in \H \;
| \; 0 \leq s \leq 1 \right\}$ connecting $p_0$ with $p_1$ is
contained in the interior of $C$. Moreover, letting $p_s = (s
z_1,s t_1)$, we have
\[
   \partial C = \{ p_0, \, p_1 \}\cup
    \bigcup_{0 < s < 1}  \partial  C \cap H_{p_s}.
\]
Note that, by strict convexity,
$ \partial C \cap H_{p_0} = \{p_0\}$ and $\partial C \cap
H_{p_1} = \{p_1\}$.

Using the property $p \in \mathrm{H}_{p'}$ if and only if $p' \in
\mathrm{H}_{p}$, we deduce that if $p \in \partial C \cap
\mathrm{H}_{p_s}$ with $0 < s < 1$, then $p_s \in \mathrm{H}_p$
and $\mathrm{H}_{p}$ is not a supporting plane of $C$, because
$p_s$ belongs to the interior of $C$. It follows that $\Sigma (C)
= \{ p_0, \, p_1 \}$.
\end{proof}

\begin{prop}\label{colx}
For all $0 < r_1 < r_2$, there is $\e > 0$, such that, for
any convex set $C \sub \H$ with $\B(p,r_1) \sub C \sub \B(p,r_2)$
for some $p \in \H$, there exist points $p_0, \, p_1 \in \Sigma (C)$
with $d(p_0,p_1) \geq \eps$.
\end{prop}
\begin{proof}
We first assume that $C$ is a strictly convex set with $C^\infty$
boundary. By Lemmas \ref{lello} and \ref{lallo}, $\Sigma(C)$
contains at least one and at most two points. After a left
translation and an isometry, we have $p_0 = 0\in\Sigma(C)$ and $C
\sub \mathrm{H}_0^+$.

There exist $r > 0$ and a smooth, strictly convex function
\[
   f:\{ z \in \R^2 \; | \; |z| < r \} \to \R,
\]
such that $f(0) = 0$, $f(z) > 0$ for $0 < |z| < r$,
and $(z,f(z)) \in \partial C$ for all $z$.

Denote by $V$ the orthogonal projection of the vector field $T$
onto $\mathrm{T} \partial C$ with respect to the left invariant
inner product $\langle \cdot \, , \cdot \rangle$, and let $W$ be
the projection of $V$ onto the $xy$-plane. Here and in
the following, when we speak of a projection onto the $xy$- or
$yt$-plane, it is understood that the projection is orthogonal
in the Euclidean sense. Thinking of $W$ as a mapping
$W:\{ z \in \R^2 \; \big| \; |z| < r \} \to \R^2$, we have
\[
   W(z)
   =
   \frac{\nabla f(z) + 2z^\perp}{1 + |\nabla f(z) + 2z^\perp|^2}
   \quad \mbox{for all $z$.}
\]
Observe that $(z,f(z)) \in \Sigma (C)$ if and only if
$\nabla f(z) + 2z^\perp  = 0$. Hence $0$ is the unique zero of
$W$ in $\{ z \in \R^2 \; \big| \; |z| < r \}$, provided that
$r > 0$ is small enough. Let
\[
   K := \left\{ z \in \R^2 \; | \; |z| \leq r/2 \right\}.
\]
Let $\partial_r f$ denote the radial derivative of $f$. By strict
convexity of $f$, the mapping $F:\partial K \to S^1$ given by
\[
   F(z)
   =
   \frac{\nabla f(z) + 2z^\perp }{|\nabla f(z) + 2z^\perp|}
\]
satisfies
\[
   F(z)\cdot \frac{z}{|z|}
   =
   \frac{\partial_r f (z)}{|\nabla f(z) + 2z^\perp|}  > 0
   \quad \mbox{on $\partial K$}.
\]
Hence $F$ and the Gauss map $G:\partial K \to S^1$ are homotopic.
It follows that
\begin{equation}
 \label{index}
   \mathrm{index}(V,0)
   =
   \mathrm{index}(W,0)
   =
   \mathrm{degree}(F,\partial K)
   =
   \mathrm{degree}(G,\partial K)
   =
   1.
\end{equation}
If $V$ had only one zero on $\partial C$, then the
Poincar\'{e}--Hopf index theorem would give
$
   \mathrm{index}(V,0) = \chi(C) = \chi(S^2) = 2,
$
where $\chi$ denotes the Euler--Poincar\'{e} characteristic. In
view of \eqref{index}, this is not possible.

Thus $\Sigma (C) = \{ p_0,p_1 \}$ consists of exactly two points.
We can assume that $p_0 = 0$. We show that $d(p_0,p_1) \geq \eps$
for some $\eps > 0$. There exists an interior point $p = (z,t) \in
C$ such that $|z| \leq \al$ and $t \geq \be$ for some $\al = \al
(r_2) < +\infty$ and $\be = \be(r_1,r_2) > 0$. Let $S = \{ (s z,s
t) \; | \; 0 \leq s \leq 1 \}$ be the line segment connecting $p_0
= 0 \in \Sigma (C)$ with $p$. As in the proof of Lemma
\ref{lallo}, we have $\Sigma (C) \cap \mathrm{H}_{(s z,s t)} =
\nul$ for all $0 < s \leq 1$. Therefore
\[
   \Sigma(C) \cap \mathrm{H}_p^- = \{ 0 \}.
\]
Moreover, $\B(0,\eps) \sub \mathrm{H}_p^-$ for some $\e > 0$
depending on $\al$ and $\be$.

Now assume that $C$ is only convex. By Lemma \ref{approx}, there
exists a sequence $\{ C_k \}_{k \in \N}$ of strictly convex sets
of class $C^\infty$ such that $\B(p,r_1/2) \sub C_k \sub
\B(p,2r_2)$ for all $k \in \N$ and $C_k \to C$ with respect to the
Hausdorff distance. From the first part of the proof, we know
that there exist sequences $\{ p_{k,0} \}_{k \in \N}$ and
$\{ p_{k,1} \}_{k \in \N}$ with $p_{k,0}, \, p_{k,1} \in \Sigma (C_k)$
and $d(p_{k,0},p_{k,1}) \geq \eps > 0$ for all $k \in \N$. Passing
to subsequences and relabelling if necessary, we can assume that
these sequences converge to limit points $p_0, \, p_1$ with the
desired properties.
\end{proof}

\begin{defn}
A line $\ell$ is horizontal if $\ell \sub \mathrm{H}_p$ for one
(equivalently: all) $p \in \ell$. A horizontal segment is a segment
of a horizontal line.
\end{defn}

\begin{thm}
   [Characteristic set]
   \label{H6}
If $C \sub \H$ is a convex open set, then $\Sigma(C)$ admits
the decomposition $\Sigma(C) = \Sigma^- \cup \Sigma^+$, where
$\Sigma^-, \, \Sigma^+$ are two closed, disjoint, horizontal
segments (possibly points). Moreover, if $0 \in \mathrm{int}(C)$,
then $\mathrm{H}_0 $ separates $\Sigma^-$ from $\Sigma^+$.
\end{thm}

\begin{proof}
Without loss of generality, we can assume $0 \in \mathrm{int}(C)$.
We define the convex sets
\[
   C^{-} = C \cap  \mathrm{H}_0^-
   \quad \mbox{and} \quad
   C^{+} = C \cap  \mathrm{H}_0^+.
\]
If  $p \in \partial C^{-} \cap \mathrm{H}_0 = \partial C^{+} \cap
\mathrm{H}_0 $ with  $p \neq 0$, then $0 \in \mathrm{H}_p \neq \mathrm{H}_0$
and consequently $\mathrm{H}_p \cap \mathrm{int}(C^\pm) \neq \nul$.
It follows that $\Sigma (C^\pm) \cap \mathrm{H}_0 = \{ 0 \}$, and,
by Corollary \ref{colx}, there exist points $p^-,p^+$ such that
$p^\pm \in \Sigma (C^\pm) \cap \mathrm{int} \left( \mathrm{H}_0^\pm \right)$.
The relative interior of the line segment connecting $0$ with $p^-$,
respectively $p^+$, is contained in the interior of $C^-$, respectively
$C^+$. Then the same reasoning as in the proof of Lemma \ref{lallo} gives
\begin{equation}
 \label{decomposition}
   \Sigma (C^-) = \{ 0 \} \cup \left( \Sigma  (C^-)
   \cap \mathrm{H}_{p^-} \right)
   \quad \mbox{and} \quad
   \Sigma (C^+) = \{ 0 \} \cup \left( \Sigma(C^+)
   \cap \mathrm{H}_{p^+} \right).
\end{equation}
Let $p \in \Sigma (C^-) \cap \mathrm{H}_{p^-}$ with $p \neq
p^{-}$. Since \eqref{decomposition} must hold with $p$ instead of
$p^-$, it follows that $\Sigma  (C^-) \cap \mathrm{H}_{p^-} $ is
contained in a horizontal segment. By compactness and convexity,
there exists a closed, bounded, horizontal segment $\Sigma^-$ such that
$\Sigma (C^-) \cap \mathrm{H}_{p^-} = \Sigma^-$. Similarly, we can show
that $\Sigma (C^+) \cap \mathrm{H}_{p^+} = \Sigma^+$ for some closed,
bounded, horizontal segment $\Sigma^+$. Neither $\Sigma^-$ nor
$\Sigma^+$ intersects $\mathrm{H}_0$.

Finally, we show that $\Sigma(C) = \Sigma^- \cup \Sigma^+$.
$\Sigma(C) \cap \mathrm{H}_0 = \nul$ implies
$\Sigma(C) \sub \Sigma^- \cup \Sigma^+$. On the other hand,
$p \in \Sigma^- \cup \Sigma^+$ implies $p \in \Sigma(C)$
(locality).
\end{proof}

We recall some definitions concerning convex functions. Let
$\O \sub \R^2$ be a convex open set and let $f:\O \to \R$ be
a convex function. The subdifferential $\partial f(q)$ of $f$
at $q \in \O$ is the (nonempty) set
\begin{equation}\label{subdifferential}
   \partial f(q)
   =
   \left\{
      q^{\ast} \in \R^2 \; | \; f(q') - f(q) \geq q^{\ast} \cdot (q' - q) \; \mbox{for all $q' \in \O$}
   \right\} .
\end{equation}
The property of belonging to the subdifferential is a local
property. It is well known that $f$ is locally Lipschitz
continuous in $\O$, that $\partial f(q) = \{ \nabla f(q) \}$ if
and only if $f$ is differentiable at $q \in \O$ with gradient
$\nabla f(q)$ and that $\nabla f \in BV_{\mathrm{loc}}(\O;\R^2)$.
Moreover, we have the following

\begin{lem}\label{compactness}
Let $\O \sub \R^2$ be a convex open set, let $f:\O \to \R$ be a
convex function and let $\{ f_k \}_{k \in \N}$ be a sequence of
convex functions $f_k:\O \to \R$ which converges to $f$ locally
uniformly in $\O$. Given a compact set $K \sub \O$, sequences $\{
q_k \}_{k \in \N}$ in $K$ and $\{ q_k^{\ast} \}_{k \in \N}$ with
$q_k^{\ast} \in \partial f_k(q_k)$ for all $k \in \N$, there exist
$q \in K$, $q^{\ast} \in \partial f(q)$ and subsequences
$\{ q_{k_l} \}_{l \in \N}$, $\{ q_{k_l}^{\ast} \}_{l \in \N}$, such
that $q_{k_l} \to q$ and $q_{k_l}^{\ast} \to q^{\ast}$.
\end{lem}

Let $C$ be a bounded, convex compact set in $\H$ with nonempty
interior. The projection $D$ of $C$ onto the $xy$-plane
\[
   D = \big\{
          z\in \R^2 \;|\; \textrm{there is $t\in\R$ such that} (z,t) \in C
       \big\}
\]
is a convex compact set in $\R^2$, and there are functions
$f,g:D\to\R$, with $f$ convex and $g$ concave, such that
$C = \big\{ (z,t)\in \H \;|\; z\in D,\, f(z)\leq t\leq g(z) \big\}$.
We have $f < g$ in $\mathrm{int}(D)$. If $z\in\partial D$,
then $f(z) \leq g(z)$, possibly $f(z) < g(z)$.

\begin{defn}
The characteristic set $\Sigma (f)$ of $f$ is the set of points $z
\in D$ such that the horizontal plane $\mathrm{H}_p$ at the point
$p = (z,f(z)) \in \R^2 \times \R$ is a supporting plane for $C$.
\end{defn}

\begin{rem}
A computation shows that
\begin{equation}\label{subdifferential1}
   z \in \Sigma(f) \cap \mathrm{int}(D)
   \iff -2z^{\perp} \in \partial f(z).
\end{equation}
\end{rem}

\begin{prop}
  [Lower bounds I]
  \label{lb1}
Let $\O \sub \mathrm{int}(D)$ be an open set and let $\{ f_{\e}
\}_{\e > 0}$ be a family of smooth convex functions in $\O$ such
that $f_{\e} \to f$ locally uniformly and $\nabla f_{\e} \to
\nabla f$ in $L_{\mathrm{loc}}^1(\O;\R^2)$ as $\e \downarrow 0$.
Then, for any compact set $K \sub \O \setminus \Sigma(f)$,
there are constants $\d > 0$ and $\e_0 > 0$ such that
\begin{equation}
   \label{Ko}
   |\nabla f_{\e}(z) + 2z^{\perp}| \geq \d
\end{equation}
for all $0 < \e < \e_0$ and all $z \in K$. Consequently,
\begin{equation}
   \label{Ok}
   | \nabla f(z) + 2z^{\perp} | \geq \d
\end{equation}
for $\L$-a.e.~$z \in K$.
\end{prop}

\begin{proof}
Suppose by contradiction that the pair $(\d,\e_0)$ does not exist.
Then, by Lemma \ref{compactness}, there are sequences $\{ f_{\e_k}
\}_{k \in \N}$, $\{ z_k \}_{k \in \N}$ and $\{ \nabla
f_{\e_k}(z_k) \}_{k \in \N}$ such that $z_k \to z \in K$ and
$\nabla f_{\e_k}(z_k) \to -2z^{\perp} \in \partial f(z)$. This
implies $z \in \Sigma(f)$ by \eqref{subdifferential1},
contradicting the assumption $K \cap \Sigma(f) = \emptyset$.
\end{proof}

Similarly, the projection $E$ of $C$ onto the $yt$-plane
\[
   E = \big\{ \z \in \R^2 \; | \;
              \textrm{there is $x\in\R$ such that} (x,\z) \in C
       \big\}
\]
is a convex compact set in $\R^2$ and there are functions $h,k:
E\to\R$, with $h$ convex and $k$ concave, such that $C = \big\{
(x,\z)\in \H \;|\;  \zeta \in E,\, h(\z)\leq x \leq k(\z) \big\}$.
We have $h < k$ in $\mathrm{int}(E)$. If $\z \in \partial E$, then
$h(\z) \leq k(\z)$, possibly $h(\z) < k(\z)$.

\begin{defn}\label{234}
The characteristic set $\Sigma (h)$ of $h$ is the set of points
$\z \in E$ such that the horizontal plane $\mathrm{H}_p$ at the
point  $p = (h(\z),\z) \in \R \times \R^2$ is a supporting plane
for $C$.
\end{defn}

\begin{rem}
A computation shows that
\begin{equation}
   \label{subdifferential2}
   \z = (y,t) \in \Sigma(h) \cap \mathrm{int}(E)
   \iff
   \mbox{$y \neq 0$ and $\big(h(\z)/y,1/(2y)\big) \in \partial h(\zeta)$}.
\end{equation}
\end{rem}

\begin{prop}
  [Lower bounds II]
  \label{lb2}
Let $\O \sub \mathrm{int}(E)$ be an open set and let $\{ h_{\e}
\}_{\e > 0}$ be a family of smooth convex functions in $\O$  such
that $h_{\e} \to h$ locally uniformly and $\nabla h_{\e} \to
\nabla h$ in $L_{\mathrm{loc}}^1(\O;\R^2)$ as $\e \downarrow 0$.
Then, for any compact set $K \sub \O \setminus \Sigma(h)$,
there are constants $\d > 0$ and $\e_0 > 0$, such that
\begin{equation}
   \label{arruh}
   |(\p_y h_{\e} - 2 h_{\e} \p_t h_{\e},1 - 2y \p_t h_{\e})| \geq \d
\end{equation}
in $K$ for all $0 < \e < \e_0$. Consequently,
\begin{equation}
   \label{hurra}
   |(\p_y h - 2 h \p_t h,1 - 2y \p_t h)| \geq \d
\end{equation}
$\L$-a.e.~in $K$.
\end{prop}
\begin{proof}
For some $\e_0 > 0$ the functions $h_{\e}$,  $0 < \e < \e_0$, are
uniformly bounded and uniformly Lipschitz continuous in some open
neighbourhood of $K$. Therefore, if $r > 0$ is sufficiently small,
we have
\[
   |1 - 2y \partial_t h_{\e}| \geq 1/2
\]
in $K \cap \{ |y| \leq r \}$ for $0 < \e < \e_0$. By \eqref{subdifferential2}
and Lemma \ref{compactness}, there exists $\d_1 > 0$ such that
\[
   |(\p_y h_{\e} - h_{\e}/y,\p_t h_{\e} - 1/(2y))| \geq \d_1
\]
in $K \cap \{ |y| \geq r \}$ for $0 < \e < \e_0$, with a possibly
smaller $\e_0$. Hence there exists $\d_2 > 0$ such that the condition
$|\p_t h_{\e} - 1/(2y)| \leq \d_2$ implies $|\p_y h_{\e} - h_{\e}/y| \geq \d_1/2$
and $|2 h_{\e}(\p_t h_{\e} - 1/(2y))| \leq \d_1/4$, whence
\[
   |\p_y h_{\e} - 2h_{\e} \p_t h_{\e}|
   \geq
   |\p_y h_{\e} - h_{\e}/y| - |2 h_{\e}(\p_t h_{\e} - 1/(2y))| \geq \d_1/4.
\]
The claim follows with $\d = \min \{ 1/2, \, \d_1/4, \, 2r\d_2 \}$.
\end{proof}

We conclude this section with the following observation.

\begin{lem}
  [Uniqueness]
  \label{uniqueness}
Let $C \sub \H$ be a convex compact set with nonempty interior.
Suppose that there exist two geodesics contained in $\p C$ and
with curvature $H > 0$ passing through a point $p\in \p C$. Then
they coincide in a neighbourhood of $p$.
\end{lem}
\begin{proof}
We can assume $p = 0$. Let $L > 0$ and let $\ga_0, \, \ga_1:[-L,L] \to \p C$
be geodesics parameterized by arc length, with curvature $H > 0$ and such that
$\ga_0(0)=\ga_1(0)=0$. If $\dot\ga_0(0)= \dot\ga_1(0)$, then $\ga_0=\ga_1$.
Assume by contradiction that $\dot\ga_0(0)\neq\dot\ga_1(0)$. Then $\HB_0$
is a supporting plane for $C$. (This is true also in the case
$\dot\ga_0(0)=-\dot\ga_1(0)$). But this is not possible, because
$\ga_0(s) \in \{ t > 0 \}$ for $s \in (0,L]$ and $\ga_0(s) \in \{
t < 0 \}$ for $s \in [-L,0)$.
\end{proof}
\section{Curvature equations for convex isoperimetric sets}

\setcounter{equation}{0}

We derive partial differential equations for certain vector fields built
from functions which parameterize the boundary of convex isoperimetric
sets. These equations state that the boundary of a convex isoperimetric set
$C\sub\H$ has constant horizontal curvature. We study two different
curvature equations: the equation for graphs of the form $t = f(x,y)$ and
the equation for graphs of the form $x=h(y,t)$. We call the number $ H =
{3\P(C)}/{4\V(C)}$ the horizontal curvature of $C$ (curvature for short).

\textbf{Graphs of the form $t = f(x,y)$.} In this subsection, we
denote the elements of $\H = \R^2 \times \R$ by $(z,t)$ with $t
\in \R$ and $z = (x,y) \in \R^2$. We write $z^\perp = (-y,x)$.
Let $C$ be a convex set in $\H$ of the form
\begin{equation}
 \label{dado}
  C = \big\{ (z,t) \in \H \;|\; z \in D, \, f(z) \leq t \leq g(z) \big\},
\end{equation}
where $D\sub\R^2$ is a convex compact set in the plane with
nonempty interior, and $-g,f:D \to \R$ are convex functions.

\begin{prop}
  [Curvature equation I]
  \label{pde1}
If $C\sub\H$ is a convex isoperimetric set with perimeter $\P(C)$ and
volume $\V(C)$, then the function $f:D \to \R$ satisfies in distributional
sense the partial differential equation
\begin{equation}
   \label{distr}
   \mathrm{div} \left( \frac{\nabla f(z)+2z^\perp} {|\nabla f(z)+2z^\perp|} \right) = H
\end{equation}
in $\mathrm{int}(D) \setminus \Sigma(f)$, with $H = 3\P(C)/4\V(C)$.
\end{prop}

\begin{proof}
By Theorem \ref{H6}, the characteristic set $\Sigma(f)\sub D$ is a
union of two (possibly empty) closed line segments. Let $\phi \in
C_c^\infty(\mathrm{int}(D)\setminus\Sigma(f))$ and for $\eps \in \R$,
consider the set
\[
      C_\eps  = \big\{ (z,t)\in \H \;|\; z\in D,\,
      f(z)+\eps \phi(z) \leq t\leq g(z) \big\}.
\]
The set $C_\eps$ is nonempty when $|\eps|$ is sufficiently small. We write
$\P(\eps) = \P(C_\eps)$, $\V(\eps) = \V(C_\eps)$ and $\Is(\eps) =
\P(\eps)^4/ \V(\eps)^3$.

If the set $C$ minimizes the isoperimetric ratio $\Is(C) =
{\P(C)^4}/ {\V(C)^3}$, then the function $\Is(\eps)$ has a minimum
at $\eps = 0$, and therefore we have the equation
\begin{equation}
 \label{cri}
    \Is'(0)= \left. \frac{\P^3}{\V^4} \big( 4 \P' \V - 3\P \V'\big)
    \right| _{\eps=0} =  0.
\end{equation}
There exists $\eps_0 > 0$ such that
\begin{equation}
 \label{vo}
       \V'(\eps)
       =
       -\int_D \phi(z) \, dz
          \quad \mbox{for $|\eps| < \eps_0$}.
\end{equation}
Moreover, denoting by
\[
   S_\eps = \big\{ (z,f(z) + \eps \phi(z)) \in \H \; | \; z \in D \big\}
\]
the graph of the function $f+\eps\phi$ and by
\[
    \nu_\eps = \frac{\big(\nabla f+\eps\nabla\phi,-1\big)\;\; }
    {\big(1+ |\nabla f+\eps\nabla\phi|^2\big)^{1/2}}
\]
the exterior unit normal to $S_\eps$ (defined $\mathcal{H}^2$-a.e.~on
$S_\eps$), from the Heisenberg Area Formula
\eqref{area} and from the standard Area Formula for graphs
of functions in Euclidean spaces, we find
\begin{equation}
 \label{pix}
 \begin{split}
  \P'(\eps) & = \frac{d}{d\eps} \int _ {S_\eps}
   \sqrt {( X\cdot \nu_\eps)^2+(Y\cdot\nu_\eps)^2} \, d \mathcal
   H^2
   \\&
   = \frac{d}{d\eps}
       \int_ D |\nabla f(z)+\eps\nabla \phi(z)+2z^\perp| \, dz.
\end{split}
\end{equation}
By Proposition \ref{lb1}, for each compact set
$K \sub \mathrm{int}(D)\setminus \Sigma(f)$, there is $\del > 0$ such that
$|\nabla f(z) + 2z^\perp| \geq \del$ for $\L$-a.e.~$z \in K$. Then
$|\nabla f(z)+\eps\nabla\phi(z)+2z^\perp|$ is essentially bounded away
from $0$ in a neighbourhood of $\mathrm{spt}(\phi)$ when $|\eps|$ is small enough
and we can interchange derivative and integral in the second line of \eqref{pix}.
At $\eps=0$ we obtain
\begin{equation}
 \label{gio}
   \P'(0)
   =\int_D  \frac{\nabla f(z)+2z^\perp}
                {|\nabla f(z)+2z^\perp|}
       \cdot  \nabla\phi(z) \, dz.
\end{equation}

From \eqref{cri}, \eqref{vo} and \eqref{gio} we find
\[
   4 \V(C)
   \int_D  \frac{\nabla f(z)+2z^\perp}{|\nabla f(z)+2z^\perp|} \cdot \nabla \phi(z) \, d z
   =
   -3\P(C) \int_D \phi(z) \, d z
\]
for arbitrary $\phi \in C_c^\infty(\mathrm{int}(D)\setminus\Sigma(f))$.
This is \eqref{distr}.
\end{proof}

\textbf{Graphs of the form $x = h(y,t)$.} In this subsection, we
denote the elements of $\H = \R \times \R^2$ by $(x,\z)$ with $x
\in \R$ and $\z = (y,t) \in \R^2$. Let $C$ be a convex set in $\H$
of the form
\begin{equation}
 \label{tratto}
  C = \big\{ (x,\z)\in \H \;|\; \z\in E,\, h(\z)\leq x\leq k(\z)
              \big\},
\end{equation}
where $E \sub \R^2$ is a convex compact set in the plane with
nonempty interior and $-k,h:E \to \R$ are convex functions.

\begin{prop}
  [Curvature equation II]
  \label{pde2}
If $C\sub\H$ is a convex isoperimetric set with perimeter $\P(C)$
and volume $\V(C)$, then the function $h:E\to\R$ satisfies in
distributional sense the partial differential equation
\begin{equation}
  \label{pox}
  (\p_y -2h \p_t)  \Big( \frac{u_1}{|u| } \Big) - 2y \p_t \Big( \frac{ u_2}{| u|} \Big) = H
\end{equation}
in $\mathrm{int}(E) \setminus \Sigma(h)$, with $H = 3\P(C)/4\V(C)$ and
\begin{equation}\label{ab}
   u = (u_1,u_2) = (h_y - 2 h h_t,1- 2y h_t).
\end{equation}
\end{prop}
\begin{proof}
By Theorem \ref{H6}, the characteristic set $\Sigma(h)$ is the
union of two (possibly empty) closed line segments. We have $u \in
BV_{\mathrm{loc}}(\mathrm{int}(E);\R^2)$ and, by Proposition
\ref{lb2}, for each compact set $K \sub (\mathrm{int}(E) \setminus
\Sigma(h))$, there exists $\del> 0$ such that $|u| \geq \del$
$\L$-a.e.~in $K$. Hence $u/|u| \in
BV_{\mathrm{loc}}(\mathrm{int}(E) \setminus \Sigma(h);\R^2)$.

Let $\phi \in C_c^\infty(\mathrm{int}(E)\setminus\Sigma(h))$. For any
$\e \in\R$ consider the set
\[
      C_{\e}  = \big\{ (x,\z)\in \H \;|\; \z\in E,\,
      h(\z) + \e \phi(\z) \leq x\leq k(\z) \big\},
\]
and let $\P(\e) = \P(C_{\e})$, $\V(\e) = \V(C_{\e})$ and
$\Is(\e) = \P(\e)^4/ \V(\e)^3$.

Denoting by $S_{\e} = \left\{ (h(\zeta) + \e \phi(\zeta),\zeta)\in
\H \;|\; \zeta \in E \right\}$ the graph of the function $h + \e \phi$
and by
\[
    \nu_{\e} = \frac{\big(-1 , \nabla h + \e \nabla\phi\big)\;\; }
    {\big(1+ |\nabla h + \e \nabla\phi|^2\big)^{1/2}}
\]
the exterior unit normal to $S_{\e}$ (defined $\mathcal H^2$-a.e.~on
$S_{\e}$), from the Heisenberg Area Formula \eqref{area} and from
the standard Area Formula, we get
\begin{equation}
 \label{pax}
 \begin{split}
  \P'(\e)
  & = \frac{d}{d\e}
  \int _ {S_\e}
     \sqrt {( X\cdot \nu_{\e})^2+(Y\cdot\nu_{\e})^2} \, d \mathcal H^2 \\
  & = \frac{d}{d\e}
  \int _ E
      \left(
         \big(1-2y (h_t + \e \phi_t)\big)^2+
         \big((h_y + \e \phi_y)- 2(h + \e \phi)
         (h_t + \e \phi_t)\big)^2
      \right)^{1/2} d\z.
\end{split}
\end{equation}
Let $\d_1 > 0$ such that $(h_y - 2hh_t)^2 + (1 - 2y h_t)^2 \geq \d_1^2$
$\L$-a.e.~in a neighbourhood $\O \sub \mathrm{int}(E) \setminus \Sigma(h)$
of $\mathrm{spt}(\phi)$. We can find $\e_0 > 0$ such that the integrand
in \eqref{pax} is larger than $\d_1/2$ $\L$-a.e.~in $\O$ whenever $|\e| < \e_0$.
Then we can differentiate under the integral sign and we obtain
\[
   \P'(0)
   =
   \int _E
   \frac{(h_y - 2hh_t)(\phi_y - 2(\phi h)_t) - 2y\phi_t (1 - 2y h_t)}
    {\sqrt{(h_y - 2hh_t)^2 + (1 - 2y h_t)^2}}
    \, d\z
\]
at $\e = 0$. A formal integration by parts yields
\begin{equation}
 \label{kia}
   \P'(0)
   =
   -\int_E
   \left(
   (\p_y - 2 h \p_t)  \Big( \frac{u_1}{|u|} \Big)
   - 2y \p_t \Big( \frac{u_2}{|u|} \Big)
   \right) \phi
   \,d\z.
\end{equation}

If the set $C$ minimizes the isoperimetric ratio $\Is(C) =
{\P(C)^4}/ {\V(C)^3}$, then, as in \eqref{cri}, we have $4 \P'(0)
\V(0) - 3\P(0) \V'(0)=0$. From this equation, from \eqref{kia} and
from
\[
   \V'(0) =  -\int_E \phi(\z) \, d\z,
\]
we get
\[
   4 \V(C) \int_E
   \left(
      (\p_y -2 h \p_t)  \Big( \frac{u_1}{|u|} \Big)
      - 2y \p_t \Big( \frac{u_2}{|u|} \Big)
   \right)\phi
   \, d\z
   =
   3 \P(C) \int _ E  \phi \, d\z
\]
for arbitrary $\phi \in C_c^\infty(\mathrm{int}(E)\setminus \Sigma(h))$. This is \eqref{pox}.
\end{proof}
\section{Improved regularity of the boundary}

\setcounter{equation}{0}

\label{regular!}

We recall some general properties of $BV$ vector fields. We refer the
reader to \cite{AFP}, chapter 3, for a detailed discussion and proofs.
Let $\O \sub \R^2$ be an open set. A vector field $u \in BV(\O;\R^2)$
has approximate limit $\bar u (q) \in \R^2$ at the point $q \in \O$ if
\begin{equation}\label{limit}
   \lim_{r \downarrow 0} \dashint_{B(q,r)} |u - \bar{u}(q)| \, d\L = 0 .
\end{equation}
Here and in the following, $B(q,r)$ denotes the open Euclidean ball
of radius $r$ centered at $q$. The approximate limit (if it exists)
is uniquely determined by \eqref{limit}. The approximate discontinuity
set of $u$ is the set $S_u$ of points $q \in \O$ at which $u$ has no
approximate limit. The set $S_u$ is a Borel set with vanishing $\L$-measure.
A point $q \in \O$ is an approximate jump point of $u$ if
there exist $u^+(q) \in \R^2$, $u^-(q) \in \R^2$ and $\nu_u(q) \in S^1$
such that $u^+(q) \neq u^-(q)$ and
\begin{equation}\label{jump1}
   \lim_{r \downarrow 0} \dashint_{B^{\pm}(q,r)} |u - u^{\pm}(q)| \, d\L = 0,
\end{equation}
where
\[
   B^{\pm}(q,r) = \big\{ q' \in B(q,r) \; | \; \pm (q' - q) \cdot \nu_u(q) > 0 \big\}.
\]
The triplet $\big( u^+(q),u^-(q),\nu_u(q) \big)$ is uniquely
determined by \eqref{jump1} up to multiplication of $\nu_u(q)$
with $-1$ and permutation of $u^+(q)$ and $u^-(q)$. The set $J_u$
of approximate jump points is an $\Haus$-rectifiable Borel subset
of $S_u$ with $\Haus \left( S_u \setminus J_u \right) = 0$.
Moreover, there exist Borel functions $u^+,u^-:J_u \to \R^2$ and
$\nu_u:J_u \to S^1$ such that \eqref{jump1} holds at every $q \in J_u$.
The precise representative $u^{\ast}:\O \to \R^2$ of $u$ is
defined as follows:
\begin{equation}\label{preciso}
   u^{\ast} (q) =
   \begin{cases}
      \bar{u}(q) & q \in \O \setminus S_u, \\
 \displaystyle
      \frac{1}{2} \big( u^+(q) + u^-(q) \big) & q \in J_u, \\
      0 & q \in S_u \setminus J_u.
   \end{cases}
\end{equation}
By the Lebesgue density theorem, $u(q) = u^{\ast}(q)$ for $\L$-a.e.~$q \in
\O$. Notice that $|u| \geq \d > 0$ $\L$-a.e.~in $\O$ implies
\begin{equation}\label{preciserep}
   |\bar{u}(q)| \geq \d \quad \mbox{for all $q \in \O \setminus S_u$}
   \quad \mbox{and} \quad
   |u^{\pm}(q)| \geq \d \quad \mbox{for all $q \in J_u$}.
\end{equation}

By the Riesz representation theorem and the Lebesgue decomposition theorem,
there exist vector-valued Radon measures $D^a u$ and $D^s u$ on $\O$ such that
$D u = D^a u + D^s u$ with $D^a u \ll \L$ and $D^s u \, \perp \,
\L$. Then we have the decomposition
\[
   D u = D^a u + D^c u + D^j u,
\]
where
\begin{equation}
 \label{AC}
  D^a u = \nabla u \, \L
\end{equation}
with $\nabla u \in L^1(\O;M^{2\times 2})$ and $M^{2 \times 2}$ is
the space of $2 \times 2$ matrices with real entries. We have $u
\in W^{1,1}(\O;\R^2)$ if and only if $D^s u$ vanishes. The
Cantor part $D^cu$ and the jump part $D^ju$ of $D^s
u$ are respectively
\[
   D^c u = D^s u  \res \O \setminus S_u
   \quad \mbox{and} \quad
   D^j u =  D^s u \res  S_u.
\]
By the representation theorem of Federer--Vol'pert and by the rank one theorem
of Alberti, the measures $D^c u$ and $D^j u$ admit the representations
\begin{equation}
 \label{mis}
\begin{split}
    D^j u & = (u^+ - u^-) \otimes \nu_u \; \Haus \res J_u
  \quad \textrm{and}
  \quad
    D^c u  = \eta \otimes \xi \; |D^c u|,
\end{split}
\end{equation}
where $|D^c u|$ is the total variation of $D^c u$ and $\eta, \,
\xi:\O \to S^1$ are suitable Borel maps. The measure $|D^c u|$ is
absolutely continuous with respect to $\Haus$.

Finally, we need the chain rule for $BV$ functions. Let $F \in
C^1(\R^2;\R^2)$ be a mapping such that $F(0) = 0$ and $|\nabla F|
\in L^{\infty}(\R^2)$. The function  $v = F \circ u$ is in
$BV(\O;\R^2)$ and $v \in L^{\infty}(\O;\R^2)$ if $u \in
L^{\infty}(\O;\R^2)$. Moreover, we have
\begin{equation}
 \label{Canto}
     D^a v  = \nabla F (u^{\ast}) \nabla u \,\L
                 \quad\textrm{and}\quad
     D^c v  = (\nabla F (u^{\ast}) \eta) \otimes \xi \; |D^c u|.
\end{equation}
A point $q \in \O$ belongs to $J_{v}$ if and only if $q \in J_{u}$
and $F(u^+(q)) \neq F(u^-(q))$. What's more
\begin{equation}
 \label{Jv}
    D^j v = \big(F(u^+) - F(u^-)\big) \otimes \nu_u \; \Haus \res
    J_u.
\end{equation}

We are now in a position to prove our regularity theorem for vector
fields with bounded variation arising from the parameterization of
the boundary of convex isoperimetric sets. Let $\O\sub\R^2$ be an
open set and let $a,b\in C(\O)$ be continuous functions.
We consider the differential operator $\LL$ defined by
\begin{equation}
 \label{DO}
   \LL u = (\p_1-a\p_2) u_1 + b \p_2 u_2.
\end{equation}
We have $\LL = \mathrm{div}$ when $a = 0$ and $b = 1$. On the other
hand, when $a = 2h$ and $b = -2y$, $\LL$ is the operator appearing
in the left hand side of \eqref{pox}.

\begin{thm}[Improved regularity]\label{Reggae}
Let $\O \sub \R^2$ be a bounded open set, let $u = (u_1,u_2)\in
BV(\O;\R^2)$ be a vector field, and let $a, b \in C(\O)$ be
continuous functions such that $b \neq 0$ in $\O$. Assume that:
\begin{enumerate}
\item
There exists $\d > 0$ such that $|u(q)| \geq \d$ for
$\L$-a.e.~$q \in \O$;
\item
$\LL u^\perp \in L^1(\O)$, where $u^\perp=(-u_2,u_1)$;
\item
$\LL \Big(\displaystyle\frac{u}{|u|}\Big) \in L^1(\O)$.
\end{enumerate}
Then $u/|u| \in W^{1,1}(\O;\R^2)$ and there exists a function
$\mu:J_u \to (0,+\infty)$, such that $u^- = \mu u^+$ on $J_u$.
\end{thm}

\begin{proof}
We claim that both the jump and the Cantor part of $D(u/|u|)$
vanish. Let $F:\R^2 \to \R^2$ be a smooth mapping such that
\begin{equation}
 \label{F}
  F(0) = 0, \quad  F(q) = \frac{q}{|q|}\; \textrm{ for } \; |q| \geq
  \frac{\d}{2}
  \quad \mbox{and} \quad
  |\nabla F |\in L^{\infty}(\R^2).
\end{equation}
Notice that, if $|q| > \d/2$, the derivative of $F$ has the form
\begin{equation}
 \label{DF}
 \nabla F (q) = \frac{1}{|q|^3} q^{\perp} \otimes q^{\perp}.
\end{equation}
By assumption i), the vector field $v = F \circ u$ is in $BV(\Omega;\R^2)$.
By \eqref{preciserep} and \eqref{Jv}, its jump set is
\[
   J_v = \big\{ q \in J_u \; | \; u^+(q)/|u^+(q)| \neq u^-(q)/|u^-(q)| \big\}
\]
and its jump part is
\begin{equation}
 \label{jippo}
   D^j v = \left( \frac{u^+}{|u^+|} - \frac{u^-}{|u^-|} \right)
           \otimes \nu_u  \; \Haus \res J_{u}.
\end{equation}
Analogously, if $D^c u = \eta \otimes \xi \; |D^c
u|$ for Borel functions $\eta,\xi:\O \to S^1$, then
\begin{equation}
 \label{cucco}
    D^c v = (\nabla F (u^*) \eta) \otimes \xi \; |D^c u|
\end{equation}
by \eqref{Canto}.

By assumption ii), the part of the measure
\[
   -\p_1 u_2 + b \p_2 u_1 + a \p_2 u_2
\]
concentrated on $J_u$ vanishes. From the formula for $D^j u$ in
\eqref{mis}, we can compute $D_k^j u_l$ for $k,l=1,2$, and we
obtain
\begin{equation}
\label{mox}
   \big( (u_2^- - u_2^+), b (u_1^+-u_1^-) + a (u_2^+-u_2^-) \big)
   \cdot \nu_u = 0
\end{equation}
$\Haus$-a.e.~on $J_u$.

By assumption iii), the part of the measure
\[
   (\p_1 - a \p_2)   (u_1/|u|) + b \p_2 (u_2/|u|)
\]
concentrated on $J_v$ vanishes. Thus, using \eqref{jippo}, we can
compute $D_k^j (u_l/|u|)$ for $k,l=1,2$, and we get
\begin{equation}\label{mix}
   \left(
   \frac{u_1^+}{|u^+|}-\frac{u_1^-}{|u^-|},
   \frac{a  u_1^- - b u_2^-}{|u^-|} - \frac{a u_1^{+}- b u_2^+}{|u^+|}
   \right)
   \cdot \nu_u = 0
\end{equation}
$\Haus$-a.e.~on $J_v$.

From \eqref{mix} and \eqref{mox}, we deduce that there exists
$\lambda \in \R$ such that the following system of equations is
satisfied $\Haus$-a.e.~on $J_v$:
\[
  \begin{cases}
   \displaystyle
   \frac{u_1^+}{|u^+|}-\frac{u_1^-}{|u^-|} = -\lambda ( u_2^+ - u_2^- ) \\
   \displaystyle
   \frac{b u_2^+ - a u_1^+}{|u^+|}
   - \frac{b u_2^- - a u_1^-}{|u^-|}
    = \lambda \big( b (u_1^+-u_1^-) + a (u_2^+-u_2^-) \big).
  \end{cases}
\]
By elementary linear algebra, using $b\neq 0$, we obtain the
equivalent systems
\[
  \begin{cases}
   \displaystyle
   \frac{u_1^+}{|u^+|}-\frac{u_1^-}{|u^-|} = -\lambda ( u_2^+ - u_2^- ) \\
   \displaystyle
    \frac{u_2^+}{|u^+|}-\frac{u_2^-}{|u^-|} = \lambda ( u_1^+ - u_1^- ),
  \end{cases}
 \quad
 \iff
 \qquad
  \left( \frac{u^+ }{|u^+ |}-\frac{u^- }{|u^- |}\right)
             \cdot \big( u^+ -u^-\big) = 0,
\]
that is $u^-\cdot u^+ = |u^-| |u^+|$, $\Haus$-a.e.~on $J_{v}$. It
follows that $u^-= \mu u^+$ $\Haus$-a.e.~on $J_{v}$ for some
$\mu :J_v\to (0,+\infty)$, and then $u^+/|u^+| = u^-/ |u^-|$ $\Haus$-a.e.~on
$J_{v}$. This proves that $\Haus(J_v) = 0$ and thus $D^j v = 0$.

Now we prove that $D^c v = 0$. By assumption ii), the Cantor part
of the measure
\[
   -\p_1 u_2 + b \p_2 u_1 + a \p_2 u_2
\]
vanishes. From $D^c u = \eta \otimes \xi \, |D^c u|$, we can
compute $D_k^c u_l$ for $k,l=1,2$, and we find
\begin{equation}\label{cat}
    (-\eta_2,b \eta_1 + a \eta_2) \cdot \xi = 0
\end{equation}
$|D^c u|$-a.e.~on $\O$. By assumption iii), the Cantor part of the measure
\[
   (\p_1 - a \p_2) (u_1/|u|) + b \p_2 (u_2 /|u|)
\]
vanishes. Thus, letting $\theta = (\theta_1,\theta_2) =
((u^*)^\perp \otimes (u^*)^\perp) \eta$, we can use
\eqref{preciserep}, \eqref{DF} and \eqref{cucco} to compute
$D_k^c (u_l/|u|)$ for $k,l=1,2$, and we find
\begin{equation}
 \label{dog}
   \big( \theta_1, b\theta_2-a\theta_1 \big) \cdot \xi = 0
\end{equation}
$|D^c u|$-a.e.~on $\O$. From \eqref{cat} and \eqref{dog}, we
deduce that there exists $\la \in \R$ such that
\[
  \begin{cases}
   \theta_1 = -\la \eta_2 \\
  b\theta_2-a\theta_1 = \la \big( b \eta_1+a \eta_2 \big)
  \end{cases}
 \quad
 \iff
 \quad \quad
   \begin{cases}
  \theta_1 = -\la \eta_2 \\
  \theta_2 = \la \eta_1.
  \end{cases}
\]
Here we used $b \neq 0$. This, in turn, is equivalent with
\[
   \theta \cdot \eta = \big( \big( (u^*)^\perp \otimes (u^*)^\perp \big) \eta \big) \cdot \eta = 0
\]
$|D^c u|$-a.e.~on $\O$. Using the identity $\big( \big( (u^*)^\perp
\otimes (u^*)^\perp \big) \eta \big) \cdot \eta = \big( (u^*)^\perp \cdot
\eta \big)^2$, we deduce that $(u^*)^\perp \cdot \eta = 0$, and thus
$\big( (u^*)^\perp \otimes (u^*)^\perp \big) \eta = 0$ as well. This
proves that the expression $(\nabla F(u^*) \eta) \otimes \xi$ in
\eqref{cucco} is equal to $0$, whence $D^c v = 0$.
\end{proof}

Theorem \ref{Reggae} has the following corollaries:

\begin{cor}\label{reg1}
Let $C \sub \H$ be a convex isoperimetric set and let $f:D \to \R$
be the function in \eqref{dado}. Then we have
\begin{equation}
   \frac{\nabla f(z) + 2z^\perp}
       {|\nabla f(z) + 2z^\perp|}
   \in W^{1,1}_{\mathrm{loc}} (\mathrm{int} (D)\setminus \Sigma(f); \R^2).
\end{equation}
\end{cor}

\begin{proof}
The vector field $ u(z) = (u_1(z),u_2(z)) = \nabla f(z) +
2z^\perp$ satisfies $\mathrm{div} \, u^{\perp} = - 4$ in
$\mathrm{int}(D)$. Moreover, by Proposition \ref{pde1}, we have
$\mathrm{div} (u/|u|) = H$ in $\mathrm{int}(D)\setminus\Sigma(f)$.
The claim follows from Theorem \ref{Reggae} with $a = 0$ and
$b = 1$.
\end{proof}

\begin{cor}
 \label{reg2}
Let $C \sub \H$ be a convex isoperimetric set and
let $h:E \to \R$ be the function in \eqref{tratto}.
Then we have
\begin{equation}
   \frac{u}{|u|}
   \in W^{1,1} _{\mathrm{loc}} (\mathrm{int}(E) \setminus (\{ y=0 \} \cup \Sigma(h));\R^2),
\end{equation}
with $u = (u_1,u_2) = (h_y  - 2 h h_t,1 - 2y h_t)$.
\end{cor}

\begin{proof}
We use Theorem \ref{Reggae} with $a = 2h$ and $b = -2y$. The vector
field $u$ is in $BV_{\mathrm{loc}}(\mathrm{int}(E);\R^2)$ and
satisfies
\[
   \LL u^\perp
   =
   (\p_y - 2h\p_t)(2yh_t - 1) - 2y\p_t(h_y - 2hh_t)
   =
   2 h_t(1  + 2y h_t)
   \in L^\infty_{\mathrm{loc}}(\mathrm{int}(E)).
\]
Moreover, by Proposition \ref{pde2}, $u/|u|$ satisfies in
distributional sense
\[
   \LL \Big(\frac{u}{|u|} \Big)
   =
   (\partial_y - 2h \partial_t)
   \Big( \frac{u_1}{|u|} \Big) - 2y \partial_t \Big( \frac{u_2}{|u|} \Big)
   =
   H
\]
in $\mathrm{int}(E) \setminus \Sigma(h)$, where $H>0$ is the curvature of $C$.
We have
\[
   u/|u| \in BV_{\mathrm{loc}}(\mathrm{int}(E) \setminus \Sigma(h))
\]
by Proposition \ref{lb2}. The claim follows from Theorem \ref{Reggae}.
\end{proof}

Observing that $2y h_y - 2h = 2y(h_y - 2 h h_t) + 2h(2y h_t - 1)$,
from Corollary \ref{reg2} we also get the following

\begin{cor}\label{reg3}
Let $C \sub \H$ be a convex isoperimetric set and let
$h:E \to \R$ be the function in \eqref{tratto}. Then
\begin{equation}
   \frac{\left( 1 - 2y h_t,2y h_y - 2h \right)}
        {\left| \left( h_y - 2 h h_t,1 - 2y h_t \right) \right|}
   \in W^{1,1}_{\mathrm{loc}} (\mathrm{int}(E) \setminus (\{ y=0 \} \cup \Sigma(h));\R^2).
\end{equation}
\end{cor}

\medskip

We conclude this section with a remark concerning smooth
regularizations of $BV$ vector fields. Let $\chi \in
C_c^\infty(\R^2)$ be a smoothing kernel with $\chi\geq 0$,
$\chi(q)=0$ if $|q| \geq 1$ and $\| \chi\|_1=1$. For $\eps > 0$,
let $\chi_\eps(q) = \eps^{-2} \chi(q/\eps)$. Let $\O$ be an
open subset of $\R^2$ and let $u \in L^1(\O;\R^2)$ be an
integrable vector field. Then, for any $q_0$ in
$\O_\eps = \{ q \in \O \; | \; \mathrm{dist}(q;\partial \O) > \eps\}$,
denote by
\begin{equation}
 \label{smutto}
   u_\eps(q_0) = \int_{\O}  \chi_\eps(q_0 - q) u(q) \, dq
\end{equation}
the standard mollification of $u$ at the point $q_0$.
If $u \in BV(\O;\R^2)$, then
\[
   \lim_{\eps \downarrow 0} u_\eps (q_0) = u^*(q_0)
   \qquad \textrm{for $\Haus$-a.e.~$q_0 \in \O$},
\]
where $u^*$ is the precise representative of $u$ defined in
\eqref{preciso}. In particular, if $q_0 \in J_u$, then
$u^*(q_0) = \frac{1}{2}(u^+(q_0)+u^{-}(q_0))$, where $u^\pm(q_0)\in\R^2$
satisfy \eqref{jump1}. We need the following variant of this property:

\begin{lem}
Let  $u \in BV(\O;\R^2)$ and let $q_1 \in \R^2$ be a fixed vector with
$|q_1| < 1$. Then we have
\begin{equation}
\label{primo}
  \lim_{\eps\downarrow 0} u_\eps (q_0 + \eps q_1) = u^*(q_0)
  \qquad \textrm{for every $q_0 \in \O \setminus S_u$}.
\end{equation}
Moreover, for all $q_0 \in J_u$, there exist numbers
$\alpha, \, \beta \in [0,1]$ (which may depend on $q_0$,
$q_1$ and $\chi$), such that $\alpha + \beta = 1$ and
\begin{equation}\label{secondo}
  \lim_{\eps\downarrow 0} u_\eps (q_0 + \eps q_1)
  =
  \alpha u^+(q_0) + \beta u^-(q_0).
\end{equation}
\end{lem}

\begin{proof}
We prove only \eqref{secondo}, the proof of \eqref{primo} being
analogous. We assume without loss of generality that the point
$q_0 \in J_u$ is $0$ and that the jump direction $\nu_u$ at
$0$ is $(0,1)$. Let $u^+ = u^+(0)$ and $u^- = u^-(0)$ be
vectors such that \eqref{jump1} holds at $0$. Let
$B^\pm_\e = \{ q \in \R^2 \; | \; |\e q_1 - q| < \e, \, q \cdot (0,\pm 1) > 0 \}$
for $\e > 0$ and let
\begin{equation}
 \label{invero}
    \alpha = \int_{B_\eps ^+} \chi_\eps(\eps q_1 - q) \, dq,
        \qquad
    \beta = \int_{B_\eps^-} \chi_\eps(\eps q_1 - q) \, dq.
\end{equation}
A change of variable shows that the definition of $\alpha$ and
$\beta$ does not depend on $\eps>0$. Using \eqref{invero}, we find
\[
\begin{split}
   u_\eps (\eps q_1)
   & =
   \int_{B_\eps^+}
   \chi_{\eps}(\eps q_1 - q) (u(q)-u^+) \, dq
   + \alpha u^+ \\
   & \quad +
   \int_{B_\eps^-}
   \chi_{\eps}(\eps q_1 - q) (u(q)-u^-)\, dq
   +\beta u^-.
\end{split}
\]
Now our claim \eqref{secondo} follows from \eqref{jump1} and from
\[
  \left|  \int_{B_\eps^\pm}
   \chi_{\eps}(\eps q_1 - q) (u(q)-u^\pm) \, dq \right|
   \leq c \, \|\chi\|_\infty \,
   \dashint_{B^\pm(0,2\eps)}|u(q)-u^\pm|\, dq.
\]
\end{proof}

\section{Stability of the curvature}

\setcounter{equation}{0}

\label{letto}

Let $\O \sub \R^2$ be a bounded open set, let $a,b\in
\mathrm{Lip}(\O)$ be Lipschitz functions and denote by $a_\e\in
C^\infty(\O_\e)$ the usual mollification of $a$. We consider the
differential operators $\LL$ and $\LL_{\e}$ defined by
\begin{equation}\label{lulu}
\begin{split}
   \LL u
   & =
   (\p_1- a \p_2) u_1 + b \p_2 u_2
   =
   \mathrm{tr}\big( \nabla u \, B\big), \\
   \LL_\e u
   & =
   (\p_1- a_\e  \p_2) u_1 + b  \p_2 u_2
   =
   \mathrm{tr}\big( \nabla u \, B_{\e} \big),
\end{split}
\end{equation}
where $u = (u_1,u_2)$,
\begin{equation}
  B = \left(\begin{array}{rc}
    1 &  0 \\
    -a  & b\
  \end{array}
 \right),
 \qquad
     B_\e = \left(\begin{array}{cc}
    1 &  0 \\
    -a_\e  & b  \
  \end{array}
 \right)
\end{equation}
and $\mathrm{tr}$ denotes the trace of the $2\times 2$ matrices
$\nabla u \, B$ and $\nabla u \, B_{\e}$.

Let $f \in \mathrm{Lip}(\O)$ be a Lipschitz function and denote by
$f_\e \in C^\infty(\O_\e)$ the mollification of $f$. We define the
vector fields
\begin{equation}\label{ueu}
   u  = (\p_1 f - a \p_2 f, b \p_2 f) + \omega
   \quad \mbox{and} \quad
   u_\e = (\p_1 f_\e- a_\e \p_2 f_\e, b \p_2 f_\e)+\omega,
\end{equation}
where $\omega\in\mathrm{Lip}(\O;\R^2)$ is a Lipschitz continuous
vector field.

\begin{thm}
  [Stability of the curvature]
  \label{mimix}
Let $\O\sub\R^2$ be a bounded open set, let
$a,b\in\mathrm{Lip}(\O)$ with $b\neq 0$ in $\O$, let
$\omega\in\mathrm{Lip}(\O;\R^2)$ and let $f\in \mathrm{Lip}(\O)$
be a function such that $\nabla f\in BV(\O;\R^2)$. Let $u\in
BV(\O;\R^2)$ and $u_\e \in BV(\O_\e;\R^2)$ be the vector fields
defined in \eqref{ueu}. Assume that:
\begin{enumerate}
\item
There exist $\d > 0$ and $\e_0 > 0$ such that $|u(q)| \geq \d$ for
$\L$-a.e.~$q \in \O$ and $|u_{\e}(q)| \geq \d$ for all $q \in
\O_{\e}$, $0 < \e < \e_0$;
\item
$\LL u^\perp \in L^1(\O)$, where $u^\perp = (-u_2,u_1)$;
\item
$\displaystyle \LL \Big(\frac{u}{|u|}\Big)\in L^1(\O)$.
\end{enumerate}
Then
\begin{equation}
   \label{ELLE}
      \lim_{\e \downarrow 0} \int_K \Big| \LL_\e
      \Big(\frac{u_{\e}}{|u_\e|}\Big)-
      \LL\Big(\frac{u}{|u|}\Big)\Big|\,
      d q = 0
\end{equation}
for any compact set $K \sub \O$.
\end{thm}

\begin{proof}
Let $F:\R^2\to\R^2$ be a smooth mapping satisfying \eqref{F}, let
$v = F \circ u$ and $v_\eps = F \circ u_\eps$. By \eqref{lulu},
assumption iii) and \eqref{Canto}, we have
\begin{equation}
 \label{grazia}
  \LL  v = \mathrm{tr}\big( (\nabla F \circ u^{\ast}) \nabla u \, B\big), \qquad
  \LL_\e  v_\e = \mathrm{tr}\big( (\nabla F \circ u_\e) \nabla u_\e \, B_\e \big),
\end{equation}
where $\nabla u \in L^1(\O;M^{2\times 2})$ is the density of the
measure $D^a u$ with respect to Lebesgue measure.  We denote by
$\nabla w \in L^1(\O;M^{2\times 2})$ the density with respect to
Lebesgue measure of the absolutely continuous part of the
distributional derivative of $w = \nabla f$. Letting  $\theta =
(-a,b)$, $\theta_\e = (-a_\e,b)$ and $w_{\e} = \nabla f_{\e}$, we
have the relations
\begin{equation}
 \label{uhha}
\begin{split}
   \nabla u  = B^T \nabla w + \Theta \qquad &\textrm{with}\quad
            \Theta = \p_2 f \, \nabla  \theta+\nabla \omega,
 \\
   \nabla u_\e  = B_\e^T \nabla w_\e + \Theta_\e
    \quad &\textrm{with}\quad
            \Theta _\e= \p_2 f_\e \, \nabla  \theta_\e+\nabla
            \omega.
\end{split}
\end{equation}

The gradient of $w_\e$ is
\[
   \nabla w_\e (q)
   =
   \int_\O \chi_\e(q-q') \, dD w
\]
for $q \in \O_{\e}$, where
\[
   D w
   =
   D^a w+ D^c w+ D^j w
   =
   \nabla w \L + D^c w + (w^+-w^-)\otimes \nu_w \Haus \res J_w,
\]
by the decomposition formulae \eqref{AC}--\eqref{mis}. Thus we can
split $\LL_\e v_\e$ in the following way:
\[
   \LL_\e v_\eps
   =
  \mathrm{tr} \left( (\nabla F \circ u_\eps) \Theta_\e B_\e \right)
   +
   \LL_\e^a v_\eps
   +
   \LL_\e^c v_\eps
   +
   \LL_\e^j v_\eps,
\]
where
\begin{equation}
 \label{divi}
\begin{split}
   \LL_\e^a v_\eps(q)
   & =
   \int_\O \chi_\eps (q-q')
   \mathrm{tr} \left( (\nabla F \circ u_\eps) B_\e^T \nabla w (q') B_\e \right)
   \, dq', \\
   \LL_\e^c v_\eps(q)
   & =
   \int_\O \chi_\eps (q-q')
   \mathrm{tr} \left( (\nabla F \circ u_\eps) B^T_\e \, d D^c w (q') B_\e \right), \\
   \LL_\e^j v_\eps(q)
   & =
   \int_{J_w} \chi_\eps(q-q')
   \left( (\nabla F \circ u_\eps) B^T_\e (w^+(q') - w^-(q')) \right)
   \cdot
   \left( \nu_w(q') B_\e \right)
   \, d \Haus(q').
\end{split}
\end{equation}
Here, unless specified otherwise, the functions depend on the
variable $q$. To get the last line of \eqref{divi}, we also used
the identities $A(\xi\otimes \eta) B_\e = (A\xi) \otimes (\eta^T
B_\e)$ and $\mathrm{tr}(\xi\otimes\eta) = \xi\cdot\eta$.

Similarly, we have the decomposition
\[
   \LL v
   =
 \mathrm{tr} \left( (\nabla F \circ u^{\ast})  \Theta B \right)
   +
   \mathrm{tr} \left( (\nabla F \circ u^{\ast}) B^{T} \nabla w B \right).
\]

We claim that
\begin{equation}
 \label{pirlow}
   \lim_{\eps \downarrow 0} \int _K
      | \LL_\e ^j v_\eps(q)| \, d q = 0.
\end{equation}
There is a constant $c > 0$ such that for $q \in K$ and for all
$0< \e < \e_0$, we have $| B_\e(q)| \leq c$. By Minkowski's
inequality and a change of variable, we find
\[
\begin{split}
   \int _K
      | \LL_\e ^j  v_\eps(q)|
      \, dq
   \leq
   c \int_{J_w} \!\! \int_K \chi_\eps(q-q')
           \left| \nabla F(u_\eps(q)) B_{\e}^{T}(q)(w^+(q')-w^-(q')) \right|
      \, dq \, d \Haus \\
   \leq
   c \, \|\chi\|_\infty
   \int_{J_w'} \int_{\{ |q| < 1 \}}
      \left| \nabla F(u_{\eps}(q' + \eps q)) B_{\e}^{T}(q'+\e q) (w^+(q')-w^-(q')) \right|
   \, dq \, d \Haus,
\end{split}
\]
where $J_w'$ is the intersection of $J_w$ with a compact set $K'
\sub \O$ which is independent of $0< \e < \e_0$ (we choose a
smaller $\eps_0$ if necessary). Since the function in the last
integral is uniformly bounded, we can invoke dominated convergence
and pass to the limit $\e \downarrow 0$ inside the integral.

If $q'\in J_w$, then, by \eqref{secondo}, we have
\[
    \lim_{\e \downarrow 0} w_\e(q' + \eps q)
    =
    \alpha w^+(q') + \beta w^-(q'),
\]
for suitable $\alpha, \, \beta \in [0,1]$ which depend on $q'$,
$q$ and $\chi$. This implies
\[
    \lim_{\eps \downarrow 0} \nabla F(u_\eps(q'+\eps q))
    =
    \nabla F(\alpha u^+(q') + \beta u^-(q')),
\]
where, according to \eqref{ueu} and \eqref{Jv}, $u^{\pm} =
(w_1^{\pm} - a w_2^{\pm}, b w_2^{\pm} )+\omega$. On the other
hand, by Theorem \ref{Reggae}, we have $u/|u|\in W^{1,1}
(\O;\R^2)$ and there exists a function $\mu:J_u \to (0,+\infty)$
such that $u^- = \mu u^+$ on $J_u$. Then the vectors $\alpha
u^+(q') + \beta u^-(q')$ and $u^+(q') - u^-(q')$ are parallel.
Now, from the identity
\[
   B^{T}(q') (w^+(q')-w^-(q')) = u^+(q')-u^-(q')
\]
and from \eqref{DF}, we get
\[
   \nabla F(\alpha u^+(q') + \beta u^-(q')) B^{T}(q') (w^+(q') - w^-(q'))
   =
   0
\]
for all $q' \in J_w = J_u$. This finishes the proof of
\eqref{pirlow}.

By dominated convergence, we find
\begin{equation}
 \label{O1}
   \lim_{\e \downarrow 0}
   \int_{K}
      \big| \mathrm{tr} \big(  (\nabla F \circ u_\e) \Theta_\e B_\e
      - (\nabla F \circ u^{\ast}) \Theta B\big) \big|
   \, dq = 0.
\end{equation}
With the short notation
\[
  \mathrm{L}_\e(q,q')
   =
   \left|
      (\nabla F \circ u_\e) B_\e^T  \nabla w(q')  B_\e - (\nabla F \circ u^{\ast}) B^T \nabla w \, B
   \right|,
\]
we have
\[
\begin{split}
   \int_K |\LL_\e^a v_\eps - \mathrm{tr}\big( (\nabla F \circ u^{\ast}) B^T \nabla w \, B\big) | \, d q
   & \leq
   c \int_{K} \int_\O \chi_\eps(q-q') \mathrm{L}_\e(q,q')  \, dq' \, dq.
\end{split}
\]
As before, a function depends on the variable $q$ unless stated
otherwise. By the triangle inequality,
\[
\begin{split}
  \mathrm{L}_\e(q,q')
   \leq &
   \left| (\nabla F \circ u_\eps) B_\e^T \nabla w (q') B_\e - (\nabla F \circ u_\e) B_\e^T \nabla w(q') B \right| \\
   & +
   \left| (\nabla F \circ u_\eps) B_\e^T \nabla w (q') B - \nabla F (u^{\ast}(q')) B^T(q') \nabla w(q') B \right| \\
   & +
   \left| \nabla F (u^{\ast}(q')) B^T(q') \nabla w (q') B - (\nabla F \circ u^{\ast}) B^T \nabla w B \right|,
\end{split}
\]
and we get
\begin{equation}
 \label{milos}
    \int_K |\LL_\e^a v_\eps -
    \mathrm{tr}\big( (\nabla F \circ u^{\ast}) B^T \nabla w \, B\big) | \, d q \leq
  c \,
  \left\{ \mathrm{I}_\eps + \mathrm{J}_\eps + \mathrm{K}_\eps \right\},
\end{equation}
where
\[
\begin{split}
   &
   \mathrm{I}_\eps
   =
   \int_\O |\nabla w(q')| \int_{K} \chi_\e(q-q') |B_\e(q)-B(q)|\, dq \, dq', \\
   &
   \mathrm{J}_\eps
   =
   \int_{K'} |\nabla w (q')| \int_{\{ |q| < 1 \}}
       \left| \nabla F(u_\eps(q'+\eps q)) B_\e^T(q'+\eps q) - \nabla  F(u^{\ast}(q')) B^T(q') \right|
   d q\, dq', \\
   &
   \mathrm{K}_\eps
   =
   \int_{K'} \int_{\{ |q| < 1 \}}
      \left| [(\nabla F \circ u^{\ast}) B^T \nabla w](q') - [(\nabla F \circ u^{\ast}) B^T \nabla w ](q' + \e q) \right|
   \, dq \, dq'.
\end{split}
\]
In the limit $\e \downarrow 0$, we have $\mathrm{I}_{\e}\to 0$ by
dominated convergence, $\mathrm{J}_\e \to 0$ by dominated
convergence and by \eqref{primo} and $\mathrm{K}_\e \to 0$ by the
continuity in the mean for $L^1$ functions.

An analogous argument shows that
\begin{equation}
 \label{pirlowa}
   \lim_{\eps \downarrow 0} \int _K |\LL_\e ^c v_\eps(q)| \, dq = 0.
\end{equation}
This concludes the proof of \eqref{ELLE}.
\end{proof}

\begin{rem}
We use Theorem \ref{mimix} in the following two situations.

1) In the coordinates $z=(x,y)$, we choose $a=0$, $b=1$ and
$\omega(z) = 2z^\perp$. The vector field in \eqref{ueu} is $u(z) =
\nabla f (z) + 2 z^\perp$ and
\[
     \LL \Big(\frac{u}{|u|}\Big) = \mathrm{div}\left( \frac{\nabla f(z) + 2 z^\perp}{|\nabla f(z) + 2
     z^\perp|}\right)
\]
is the curvature operator in \eqref{distr}.

2) In the coordinates $\z = (y,t)$, we have a function $h\in
\mathrm{Lip}(\O)$ and we choose $a= 2h$, $b= -2y$, and $\omega =
(0,1)$. The vector field in \eqref{ueu} is $u = (u_1,u_2) = (h_y -
2 h h_t,1- 2y h_t)$ -- this is \eqref{ab} -- and
\[
   \LL  \Big(\frac{u}{|u|}\Big) =
 (\p_y -2h \p_t)  \Big( \frac{u_1}{|u| } \Big) - 2y \p_t \Big( \frac{ u_2}{| u|} \Big)
\]
is the curvature operator in \eqref{pox}.
\end{rem}


\section{Foliation by geodesics}

\label{Lag}

Let $u \in BV(\R^2;\R^2) \cap L^{\infty}(\R^2;\R^2)$ and, for $\rho > 0$
and $q \in \R^2$, define
\[
   C_q([-\rho,\rho];\R^2) =
   \left\{
      \ga \in C([-\rho,\rho];\R^2) \; \Big| \;
      \ga(s) = q + \int_0^s u(\ga(\s)) \, d\s, \; s \in [-\rho,\rho]
   \right\}.
\]
Both $C([-\rho,\rho];\R^2)$ and $C_q([-\rho,\rho];\R^2)$ are endowed with
the sup norm. Let $K\sub\R^2$ be a compact set. An $\L$-measurable map
$\Phi: K \to C([-\rho,\rho];\R^2)$ is a Lagrangian flow starting from $K$
relative to $u$ if $\Phi(q) \in C_q([-\rho,\rho];\R^2)$ for $\L$-a.e.~$q\in
K$. By abuse of notation, we identify the map $\Phi$ with the map $\Phi :
K  \times[-\r,\r]\to \R^2$, $\Phi(q,s) = \Phi(q)(s)$. The flow is said to
be regular if there exists a constant $m\geq 1$ such that
\begin{equation}
 \label{regular}
  \frac{1}{m} \L(A) \leq \L(\Phi (A,s))
   \leq m \L(A)
\end{equation}
for all $\L$-measurable sets $A \sub K$ and for all $s \in [-\rho,\rho]$.

\begin{thm}[Ambrosio]
  \label{flow}
Assume that:
\begin{enumerate}
\item
   $u \in BV(\R^2;\R^2) \cap L^{\infty}(\R^2;\R^2)$;
\item
   $\mathrm{div} \, u \in L^{\infty}(\R^2)$.
\end{enumerate}
Then, for any compact set $K\sub\R^2$ and for all $\rho>0$, there exists a
unique regular Lagrangian flow $\Phi:K\times[-\rho,\rho]\to\R^2$ starting
from $K$ relative to $u$, unique in the sense that any other regular
Lagrangian flow for $u$ starting from $K$ coincides with $\Phi$ $\L$-a.e.
Moreover:
\begin{enumerate}
\item[1)] The flow is stable, in the sense that, if $\{ u_{\eps} \}_{0 <
\eps < \eps_0}$ is a family of smooth vector fields which is uniformly
bounded in $BV(\R^2;\R^2) \cap L^{\infty}(\R^2;\R^2)$ and such that the
$u_{\e}$ have uniformly bounded divergence, then $u_{\e} \to u$ in
$L^1(\R^2;\R^2)$ implies
\begin{equation}
 \label{stability}
   \lim_{\e \downarrow 0}
      \int_K \sup_{s \in [-\r,\r]} |\Phi_{\e}(q,s) - \Phi(q,s)|
   \, dq =0,
\end{equation}
where $\Phi_\e$ denotes the flow of $u_\e$. In particular, the mapping $
(q,s) \mapsto \Phi(q,s)$ is $\L \times \mathcal{L}^1$-measurable.

\item[2)] The flow satisfies the semigroup property
\begin{equation}\label{semigroup}
   \Phi(\Phi(q,s),\s) = \Phi(q,s + \s)
\end{equation}
for all $s, \s \in [-\rho,\rho]$ with $s + \s \in [-\rho,\rho]$.
\end{enumerate}

\end{thm}

The existence statement is Theorem 6.2, the uniqueness statement is Theorem
6.4, the stability statement is Theorem 6.6, the semigroup property is
Remark 6.7 in \cite{A}. Ambrosio's theory holds more generally for non
autonomous vector fields in any space dimension.

\begin{rem}\label{compositions with the flow}
\hfill\newline 1) If $\O \sub \R^2$ is an open set, $K \sub \O$ is compact
and $u \in BV(\O;\R^2) \cap L^{\infty}(\O;\R^2)$ with $\mathrm{div} \, u \in
L^{\infty}(\O)$, then there is  a vector field in $BV(\R^2;\R^2) \cap
L^{\infty}(\R^2;\R^2)$ with distributional divergence in
$L^{\infty}(\R^2)$, which coincides with $u$ in a neighbourhood of $K$. This
vector field has a regular Lagrangian flow $\Phi: K\times [-\r,\r] \to
\R^2$. If $\r > 0$ is sufficiently small, then, for almost every $q \in K$,
$s \mapsto \Phi(q,s)$ is an integral curve of $u$ and $\Phi(q,s) \in \O$
for all $s \in [-\r,\r]$.

2) If $N \sub \R^2$ is negligible, then $\Phi^{-1}(N) \sub K \times [-\r,\r]$
is also negligible. This follows from \eqref{regular} by a Fubini-type argument.

3) If $F, \, G:\R^2 \to \R^k$ are $\L$-measurable mappings such that
$F = G$ $\L$-a.e., then, for $\L$-a.e.~$q \in K$,
$F(\Phi(q,s)) = G(\Phi(q,s))$ for a.e.~$s \in [-\r,\r]$.

4) If $F:\R^2 \to \R$ is locally Lipschitz continuous, then the chain rule
holds along integral curves of $u$, i.e.
\begin{equation}\label{chain}
   \frac{d}{ds}F(\Phi(q,s))
   =
   \nabla F(\Phi(q,s)) \, \frac{d}{ds} \Phi(q,s)
\end{equation}
for $\L$-a.e.~$q \in K$ and for a.e.~$s \in [-\r,\r]$.
\end{rem}

We are interested in the flow associated with a certain horizontal
vector field tangent to the boundary of a convex isoperimetric set
$C \sub \H$. Let $f:D \to \R$ be the convex function in
\eqref{dado}. We consider the vector fields $v(z) = 2z - \nabla
f^\perp(z)$ and $u(z) = v^{\perp}(z) = \nabla f(z) + 2z^{\perp}$
in $BV_{\mathrm{loc}}(\mathrm{int}(D);\R^2) \cap
L^\infty_{\mathrm{loc}}(\mathrm{int}(D);\R^2)$. We have
\begin{equation}\label{coop}
   \mathrm{div} \, v
   =
   4 + f_{yx} - f_{xy}
   =
   4.
\end{equation}
The vector field $v$ is the projection onto the $xy$-plane of the vector
field
\[
   (x,y) \mapsto
   (f_y + 2x)\partial_x + (2y-f_x)\partial_y + (2xf_x + 2y f_y)\partial_t,
\]
which is both horizontal and tangent to the graph of $f$ at
$\mathcal{H}^2$-a.e.~point.

Equation \eqref{distr} completely determines the geometry of the integral
curves of $v$.

\begin{thm}
  [Foliation by circles I]
  \label{effe}

Let $C \sub \H$ be a convex isoperimetric set with curvature $H >
0$ and $f:D \to \R$ be the function in \eqref{dado}. Let $K \sub
\mathrm{int}(D) \setminus \Sigma(f)$ be a compact set, $\O \Subset
\mathrm{int}(D) \setminus \Sigma(f)$ be an open neighbourhood of
$K$ and $\r > 0$ be sufficiently small. Denote by $\Phi:K \times
[-\rho,\rho] \to \O$ the regular Lagrangian flow relative to the
vector field $v(z) = 2z- \nabla f^\perp(z)$ starting from $K$.
Then, for $\L$-a.e.~$z \in K$, the integral curve $s \mapsto
\Phi(z,s)$ is an arc of circle with radius $1/H$ oriented
clockwise.
\end{thm}

In the case of graphs of the type $x = h(y,t)$, the result is similar. Let
$h:E \to \R$ be the function in \eqref{tratto}. The vector field
\begin{equation}
 \label{uva}
   v(\z) = \big(1 - 2y h_t,2y h_y - 2h\big),
   \quad
   \z = (y,t) \in E ,
\end{equation}
is in $BV_{\mathrm{loc}}(\mathrm{int}(E);\R^2)\cap
L^\infty_{\mathrm{loc}}(\mathrm{int}(E);\R^2)$, and
\begin{equation}
 \label{migros}
   \mathrm{div}\, v = - 4 h_t -2y h_{ty} +2y h_{yt}= -4 h_t
   \in L^\infty_{\mathrm{loc}}(\mathrm{int}(E)).
\end{equation}
The vector field $v$ in \eqref{uva} is the projection onto the $yt$-plane
of the vector field
\begin{equation}
 \label{cosca}
   (y,t)
   \mapsto
     (h_y - 2 h h_t) \partial_x
   + (1 - 2y h_t) \partial_y
   + (2y h_y - 2h) \partial_t,
\end{equation}
which is both horizontal and tangent to the graph of $h$ at
$\mathcal{H}^2$-a.e.~point.

\begin{thm}
  [Foliation by circles II]
 \label{acca}
Let $C \sub \H$ be a convex isoperimetric set with curvature $H > 0$ and
$h:E\to \R$ be the function in \eqref{tratto}. Let
$K \sub \mathrm{int}(E) \setminus \big( \{ y = 0 \} \cup \Sigma(h)\big)$
be a compact set, $\O \Subset \mathrm{int}(E) \setminus \big( \{ y = 0 \} \cup \Sigma(h)\big)$
be an open neighbourhood of $K$ and $\r > 0$ be sufficiently small. Denote by
$\Phi: K\times[-\rho,\rho] \to \O$ the regular Lagrangian flow relative to the
vector field $v$ in \eqref{uva} starting from $K$. Then, for $\L$-a.e.~$\z \in K$,
the projection onto the $xy$-plane of the curve
\begin{equation}
 \label{migrol}
    s \mapsto ( h(\Phi(\z,s)),\Phi(\z,s)) \in \H,
    \qquad
    s \in [-\rho,\rho],
\end{equation}
is an arc of circle with radius $1/H$ oriented clockwise.
\end{thm}

We give two different proofs of Theorems \ref{effe} and \ref{acca} in
Sections \ref{xillo} and \ref{Sobolev}, respectively. The proof in Section
\ref{xillo} is based on Theorem \ref{mimix} and on the stability property
\eqref{stability} of the flow, while the proof in Section \ref{Sobolev}
relies on a reparameterization argument involving \eqref{regular}. These
theorems have the following geometric interpretation.

\begin{cor}[Foliation by geodesics I]
 \label{geodesics1}
Under the assumptions of Theorem \ref{effe}, for all $z \in K$, there is a
geodesic $\ga_z:[-\r,\r]\to\p C$ with curvature $H$, such that $\ga_z(0) =
(z,f(z))$. Moreover, the length of $\ga_z$ is bounded from below by a
positive constant depending on $K$.
\end{cor}
\begin{proof}
The horizontal lift of the plane curve $s \mapsto \Phi(z,s)$ given
by Theorem \ref{effe}, i.e.~the curve
\[
   \ga_z(s) = (\Phi(z,s),f(\Phi(z,s))),\quad s\in[-\r,\r],
\]
is a geodesic with curvature $H$ and with length bounded from
below by a positive constant depending on $K$. By Theorem
\ref{effe}, this curves exist for $\L$-a.e.~$z\in K$. The stated
result now follows from a density-compactness argument.
\end{proof}

Analogously, we have:

\begin{cor}[Foliation by geodesics II]
  \label{geodesics2}
Under the assumptions of Theorem \ref{acca}, for all $\z \in K$, there is
a geodesic $\ga_\z:[-\r,\r]\to\p C$ with curvature $H$, such that
$\ga_\z(0) = (h(\z),\z)$. Moreover, the length of $\ga_\z$ is bounded
from below by a positive constant depending on $K$.
\end{cor}
\section{Approximation argument}

\setcounter{equation}{0}

\label{xillo}

We prove Theorems \ref{effe} and \ref{acca} using the stability property
\eqref{stability} of the flow and the stability property \eqref{ELLE} of
the curvature with respect to smooth approximations.

\begin{proof}[Proof of Theorem \ref{effe}]
For $\e > 0$ sufficiently small, the mollification $f_{\e}$ of $f$
is defined in some open neighbourhood of $\overline{\O}$. Let
$v_{\e} = 2z-\nabla f_{\e}^{\perp}$. From \eqref{coop}, it follows
that
\begin{equation}
  \label{quaqua}
   \mathrm{div} \, v_{\e} = \mathrm{div} \, v = 4 \quad \mbox{in $\O$}.
\end{equation}
In particular, $\mathrm{div} \, v_{\e}, \, \mathrm{div} \, v \in L^1(\O)$.
Note that, by Proposition \ref{lb1}, there are constants $\d > 0$
and $\e_0 > 0$, such that $|\nabla f_{\e} (z)+ 2 z^{\perp}| \geq \d$ for
all $z \in \O$ and $0 < \e < \e_0$, and $|\nabla f(z)+ 2 z^{\perp}| \geq \d$
for $\L$-a.e.~$z \in \O$.  Let $F \in C^{\infty}(\R^2;\R^2)$ be a
mapping which satisfies \eqref{F}. We define $w_{\e} = F \circ v_{\e}$
in $\O $  and $w = F \circ v$ in $\mathrm{int}(D)$. Then
$w_{\e} = v_{\e}/|v_{\e}|$ and $w = v/|v|$ in $\O$. Moreover,
\begin{equation}
   \mathrm{div}\, w^{\perp} = H
\end{equation}
in the weak sense, by Proposition \ref{pde1} and Corollary \ref{reg1}.

Denote by $\Phi_{\e} :K  \times [-\r,\r] \to \R^2$ the flow of $v_{\e}$. We
have $\Phi_{\e}(K \times [-\r,\r]) \sub \O$ when $\r$ is sufficiently small
and $0 < \e < \e_0$. Denote by $J_\e(s) = J\Phi_\e(\cdot,s)$ the Jacobian
in the $z$-variable of the flow. $J_\e$ satisfies the differential equation
$\dot J_\e  = \mathrm{div}\, v_\e J_\e$, where the dot denotes the
derivative with respect to the variable $s$. We find
\begin{equation}
 \label{Jace}
   J\Phi_\e(z,s)  
   = e^{4s}
   \quad \mbox{for all $s \in [-\r,\r]$}.
\end{equation}

We define the curvature function $H_\eps:K \times [-\r,\r] \to \R$ to be
\begin{equation}
 \label{H_e}
         H_\e(z,s) =  \mathrm{div}\,
         w_\e^\perp (\Phi_\e(z,s))
         \quad \mbox{for all $z \in K$ and $s \in [-\r,\r]$}.
\end{equation}

\noindent For $z \in K$, the integral curve
$s \mapsto \Phi_\eps(z,s)$ is of class $C^\infty$. The unit tangent
vector $\GG_\e(z,s) = F(\dot\Phi_\e(z,s)) = w_\e(\Phi_\e(z,s))$
satisfies the equation
\begin{equation}\label{dt}
   \dot  \GG_\e(z,s)
   =
   - H_\e(z,s)  \la_\e(z,s)  \GG_{\e}^{\perp}(z,s)
   \quad \mbox{for all $z \in K$ and $s \in [-\r,\r]$},
\end{equation}
where $\la_\e(z,s) = |v_\e(\Phi_\e(z,s))|$. From \eqref{dt}, we
obtain
\begin{equation}
  \label{dtdt}
  \int_K \psi  \int_{-\r }^\r \Big\{
  \GG_\e \cdot\dot\phi -H_\e \la_\e
  \GG_{\e}^{\perp} \cdot \phi \Big\}
  \, ds\, dz
  = 0
\end{equation}
for arbitrary test functions $\phi \in C_c^1((-\r,\r);\R^2)$ and
$\psi \in C_c^1(K)$. Our aim is to prove that
\begin{equation}
 \label{prima}
   \int_K \psi
   \int_{-\r}^\r
   \Big\{ \GG  \cdot \dot \phi
   -H \la \, \GG^{\perp} \cdot \phi
   \Big\}\, d s\, d z  = 0,
\end{equation}
with $\GG(z,s) = w(\Phi(z,s))$ and $\la(z,s) = |v(\Phi(z,s))|$.

We need some preliminary observations. By the stability property
\eqref{stability} of the Lagrangian flow, there is a sequence $\{
\e_k \}_{k \in \N}$ such that $\e_k \downarrow 0$ and $\Phi_{\e_k}
\to \Phi$ pointwise almost everywhere in $K \times [-\r,\r]$. We
claim that
\begin{equation}
  \label{topo}
   \lim_{k \to \infty}
   \int_K \int_{-\r}^{\r}
   \left| w_{\e_k}\circ \Phi_{\e_k} - w \circ \Phi \right| \, d s \, d z
   = 0.
\end{equation}
Indeed, by \eqref{Jace}, we have
\begin{equation}
   \label{gatto}
   \int_K \left| w_\e\circ \Phi_\e  - w \circ \Phi_\e\right| \, dz
   \leq e^{-4s} \int_{\O} | w_\e - w | \, dz,
\end{equation}
and the right hand side of \eqref{gatto} goes to zero. Moreover,
we have
\begin{equation}
  \label{quatto}
   \lim_{k \to \infty}
   \int_K \int_{-\r}^\r
   | w \circ \Phi_{\e_k} - w \circ \Phi | \, d s \, d z
   = 0.
\end{equation}
Indeed, let $\eta > 0$ be a positive number and choose a vector field
$\tilde{w} \in C_c(\O;\R^2)$ such that $\| w - \tilde{w} \|_{L^1(\O;\R^2)} \leq \eta$.
Denote by $\| \cdot \|$ the norm in $L^1(K \times [-\r,\r];\R^2)$. By \eqref{regular}
and \eqref{Jace}, we have
\[
\begin{split}
   \| w \circ \Phi_\e - w \circ \Phi \|
   \leq &
   \| w \circ \Phi_\e - \tilde{w} \circ \Phi_\e \|
   + \| \tilde{w}  \circ \Phi_\e - \tilde{w} \circ \Phi \|
   + \| \tilde{w} \circ \Phi - w \circ \Phi \| \\
   \leq &
   c \eta + \| \tilde{w} \circ \Phi_\e - \tilde{w} \circ \Phi \|,
\end{split}
\]
where $c > 0$ does not depend on $\e$. This finishes the proof of
\eqref{topo}. From \eqref{topo}, it follows that
\begin{equation}
 \label{quarto}
 \begin{split}
   &
   \lim_{k \to \infty}
   \int_K \psi  \int_{-\r}^\r
   H \la \, \GG_{\e_k}^{\perp} \cdot \phi  \, d s \, d z
   =
   \int_K \psi  \int_{-\r}^\r
   H \la \, \GG^{\perp} \cdot \phi
      \, ds \, dz,
 \\
 &
   \lim_{k\to\infty}
   \int_K
   \psi \int_{-\r}^{\r} \, \GG_{\e_k}  \cdot \dot{\phi}
   \, ds \, dz
   =
   \int_K   \psi\int_{-\r}^{\r}
   \, \GG \cdot \dot{\phi}
   \, ds \, dz.
\end{split}
\end{equation}

Finally, performing the change of variable $z \mapsto
\Phi_{\e}(z,-s)$ with area factor given by \eqref{Jace}, we obtain
\begin{equation}
  \label{kronprinz}
   \int_K |(\mathrm{div}\, w^\perp_\e) \circ \Phi_\e -H| \, dz
   \leq
   e^{-4s} \int_{\O} |\mathrm{div}\, w^\perp _\e -H| \, dz,
\end{equation}
for any $s \in [-\r,\r]$. Due to \eqref{ELLE} in Theorem
\ref{mimix} with $a=0$, $b=1$ and $\omega=2z^\perp$, the right
hand side of \eqref{kronprinz} converges to $0$ as $\e \downarrow
0$ (we have $\LL =\mathrm{div}$). From \eqref{kronprinz},
we deduce
\begin{equation}
 \label{koenig}
    \lim_{k\to\infty}
    \int_K \psi \int_{-\r}^{\r}
 \left( H_{\e_k} \la_{\e_k} - H \la  \right) \,
    \GG_{\e_k}^{\perp} \cdot \phi
    \, ds \, dz = 0.
\end{equation}

Now, \eqref{prima} follows from \eqref{quarto} and \eqref{koenig}.

Since $\psi\in C_c^1(K)$ is arbitrary in \eqref{prima}, we deduce
that the equation
\begin{equation}\label{cpcp}
   \int_{-\r}^\r
   \Big\{  \GG (z,s) \cdot \dot \phi(s)
   - H \la(z,s) \GG^{\perp}(z,s) \cdot \phi(s)
   \Big\}
   \, d s = 0,
\end{equation}
holds for $\L$-a.e.~$z \in K$ and for all $\phi \in
C_c^1((-\r,\r);\R^2)$. Now, by a standard argument, we conclude
that for $\L$-a.e.~$z \in K$ we have
\[
      \dot\Gamma(z,s) = -  H \la(z,s) \GG^{\perp}(z,s)
\]
for a.e.~$s\in[-\r,\r]$, where $\la(z,t) =|v(\Phi(z,t))|=
|\dot\Phi(z,t)|$ is the length factor of the curve $s\mapsto
\Phi(z,s)$. Then $\Phi(z,\cdot)$ parameterizes an arc of circle
with curvature $H$.

\end{proof}

\begin{proof}[Proof of Theorem \ref{acca}]
For $\e > 0$ sufficiently small, the mollification $h_{\e}$ of $h$ is defined
in some open neighbourhood of $\overline{\O}$. We consider the vector fields
\[
   v_{\e} = (1 - 2y \p_t h_{\e},2y \p_y h_{\e} - 2h_{\e})
   \quad \mbox{and} \quad
   u_{\e} = (\p_y h_{\e} - 2h_{\e} \p_t h_{\e},1 - 2y \p_t h_{\e}).
\]
Analogously, in $\mathrm{int}(E)$, we can define
\[
   v = (1 - 2y h_t,2y h_y - 2h)
   \quad \mbox{and} \quad
   u = (h_y - 2h h_t,1 - 2y h_t).
\]
By Proposition \ref{lb2}, there are constants $\d > 0$ and $\e_0 > 0$, such
that $|u_{\e}| \geq \d$ in $\O$ for $0 < \e < \e_0$, and $|u| \geq \d$
$\L$-a.e.~in $\O$. Let $F \in C^{\infty}(\R^2;\R^2)$ be a mapping which
satisfies \eqref{F}. We can define the vector fields $w_{\e} = F \circ u_{\e}$
in $\O$ and $w = F \circ u$ in $\mathrm{int}(E)$. Then $w_{\e} = u_{\e}/|u_{\e}|$
and $w = u/|u|$ in $\O$.

Let $\LL$ be the differential operator in \eqref{DO}, with $a =
2h$ and $b = -2y$. By Proposition \ref{pde2} and Corollary
\ref{reg2}, the vector field $u$ satisfies
\[
   \LL \Big( \frac{u}{|u|} \Big) = H
   \quad \mbox{in $\O$},
\]
in the weak sense. A direct computation shows that
\[
   \LL u^{\perp} = 2 h_t(1 + 2y h_t) \in L^\infty(\O).
\]
Moreover, from \eqref{migros}, we get
\[
   \mathrm{div} \, v_\e = - 4 \p_t h_\e
   \quad \mbox{in $\O$},
\]
and there is a constant $c > 0$ such that $|\mathrm{div} \, v_{\e}| \leq c$
in $\O$ for all $0 < \e < \e_0$. Let us denote by $\Phi_{\e}: K \times
[-\r,\r] \to \R^2$ the flow of $v_{\e}$. As in the proof of Theorem
\ref{effe}, $\r$ can be chosen in such a way that $\Phi_{\e}(K \times
[-\r,\r]) \sub \O$ for $0 < \e < \e_0$, and the Jacobian of $\Phi_{\e}$
satisfies
\begin{equation}
 \label{Jace2}
    c_1 \leq J \Phi_{\e}(\z,s) \leq c_2
    \quad \mbox{for all $\z \in K$ and $s \in [-\r,\r]$},
\end{equation}
with constants $0 < c_1 \leq c_2$ which do not depend on $\e$.

Let $s \mapsto \GG_\e(\z,s) = w_\e(\Phi_\e(\z,s))$ be the unit
tangent vector of the projection onto the $xy$-plane of the curve
$s\mapsto \big(h_\e(\Phi_\e(\z,s)), \Phi_\e(\z,s)\big)\in\H$. This
tangent vector satisfies the equation
\begin{equation}
 \label{dtds}
 \dot  \GG_\e(\z,s)
   = - H_\e(\z,s) \la_\e (\z,s) \GG_\e ^\perp(\z,s)
   \quad \mbox{for all $\z \in K$ and $s \in [-\r,\r]$},
\end{equation}
where $H_\e:K\times[-\r,\r]\to \R$ is the curvature function
\begin{equation}
 \label{kux}
   H_\e(\z,s) =  \LL\Big( \frac{u_\e}{|u_\e|}\Big)(\Phi_\e(\z,s)),
\end{equation}
and $\la_\e(\z,s) = |u_\e(\Phi_\e(\z,s))|$. The proof of formula
\eqref{dtds} is postponed to Section \ref{Sobolev}, where it is
proved under weaker regularity assumptions (see the proof just
before formula \eqref{dirx}).

As in \eqref{dtdt}, we arrive at the integral equation
\begin{equation}
 \label{dtdtds}
  \int_K \psi  \int_{-\r }^\r
     \Big\{
            \GG_\e \cdot \dot \phi - H_\e \la_\e \, \GG_\e ^\perp \cdot \phi
     \Big\}
   \, ds \, d\z = 0,
\end{equation}
with $\phi\in C^1_c((-\r,\r);\R^2)$ and $\psi\in C_c^1(K)$. Now
the proof runs along the same lines as the proof of Theorem
\ref{effe} above, and we can omit the details. Theorem \ref{mimix}
is used with $a=2h$, $b=-2y$ and $\omega=(0,1)$.
\end{proof}


\section{Reparameterization argument}

\setcounter{equation}{0}

\label{Sobolev}
Let $\O\sub\R^2$ be an open set and let $v$ be a vector field such that
\begin{equation}
 \label{hippy}
      v\in BV(\O;\R^2) \cap L^{\infty}(\O;\R^2)\quad
      \textrm{ and }
      \quad
      \mathrm{div} \, v \in L^{\infty}(\O).
\end{equation}
The function $v$ is defined pointwise, i.e.~we choose a
representative in the equivalence class of $v$. Our results hold
independently of this choice. However, there is an exceptional set
of points which may a priori depend on the representative.

Given a compact set $K \sub \O$ and a sufficiently small $\r > 0$,
there exists a unique regular Lagrangian flow $\Phi: K \times[-\r,\r] \to \O$
starting from $K$ relative to $v$ (Theorem \ref{flow}). Let $\la:\O\to\R$
be a measurable function such that
\begin{equation}
 \label{ciuno}
   0 < c_1 \leq \la \leq c_2 \quad \textrm{$\L$-a.e.~in $\O$}.
\end{equation}
Then, for $\L$-a.e.~$q \in K$, the curve
$s \mapsto \lambda(\Phi(q,s))$ is measurable and
\[
   c_1 \leq \la(\Phi(q,s)) \leq c_2
   \quad \mbox{for a.e.~$s \in [-\rho,\rho]$}.
\]

Thus, for $\L$-a.e.~$q \in K$, the change of parameter $\sigma_q:[-\r,\r] \to
[\sigma_q(-\r),\sigma_q(\r)]$ defined by
\[
   \sigma_q(s) = \int_0^s \la(\Phi(q,\xi)) \, d\xi
\]
is bi-Lipschitz, strictly increasing and admits therefore a
bi-Lipschitz, strictly increasing inverse
$\tau_q:[\sigma_q(-\r),\sigma_q(\r)] \to [-\r,\r]$,
which satisfies
\begin{equation}
 \label{dif}
   \dot{\tau}_q(s)
   =
   \frac{1}{\dot{\sigma}_q(\tau_q(s))}
   =
   \frac{1}{\la (\Phi(q,\tau_q(s)))}
   \quad \mbox{for a.e.~$s \in [\sigma_q(-\r),\sigma_q(\r)]$}.
\end{equation}
Consequently, for $\L$-a.e.~$q \in K$, the curve
$\ga_q:[\sigma_q(-\r),\sigma_q(\r)] \to \O$ defined by
\begin{equation}
 \label{gaia}
   \ga_q(s) = \Phi(q,\tau_q(s))
\end{equation}
is absolutely continuous and satisfies
\begin{equation}\label{gioia}
   \dot{\ga}_q(s) = \frac{v(\ga_q(s))}{\la(\ga_q(s))}
   \quad \mbox{for a.e.~$s \in [\sigma_q(-\r),\sigma_q(\r)]$},
\end{equation}
i.e.~$\ga_q$ is an integral curve of the vector field $v/\la$.

\begin{thm}[Chain rule for integral curves]

 \label{regularity}

Let $\O\sub\R^2$ be an open set, let $K\sub\O$ be a compact subset, let $w
\in W^{1,1}(\O;\R^2)$ and let $\r > 0$ be sufficiently small. Then, for
$\L$-a.e.~$q \in K$, the curve $w \circ \ga_q$ with $\ga_q$ defined in
\eqref{gaia} belongs to $W^{1,1}((\sigma_q(-\r),\sigma_q(\r));\R^2)$, and
its weak derivative is $(\nabla w \circ \gamma_q) \dot{\ga}_q$.
\end{thm}

\begin{proof}
Let $\{ w_k \}_{k \in \N}$ be a sequence of smooth vector fields $w_k\in
W^{1,1}(\O;\R^2)$ converging to $w$ in $W^{1,1}(\O;\R^2)$. Notice that the
map $(q,s) \mapsto w(\Phi(q,s))$ is measurable and integrable. By Fubini's
theorem and by the bounded volume distortion property \eqref{regular} of
the flow, we have
\[
\begin{split}
   \int_{K} \int_{-\r}^{\r}
      |w_k(\Phi(q,s)) - w(\Phi(q,s))|
   \, d s \, d q  &
   \leq
   m \,
   \int_{-\r}^{\r} \int_{\Phi (K,s)}
      |w_k(q) - w(q)|
   \, d q \, d s,
\end{split}
\]
and
\[
\begin{split}
   \int_{K} \int_{-\r}^{\r}
      \big|\big( \nabla w_k (\Phi(q,s)) - & \nabla w (\Phi(q,s))\big)
       v(\Phi(q,s)) \big|
   \, d s \, d q
     \leq
     \\
     & \leq
    \| v \|_{\infty}
   \int_{-\r}^{\r} \int_{K}
      |  \nabla w_k (\Phi(q,s)) - \nabla w (\Phi(q,s)) |
   \, d q \, d s \\
   & \leq
 m  \| v \|_{\infty}
   \int_{-\r}^{\r} \int_{\Phi (K,s)}
      | \nabla w_k(q) - \nabla w(q) |
   \, d q \, d s.
\end{split}
\]
Notice that the curve $s \mapsto w_k(\Phi(q,s))$ is Lipschitz, for all $k
\in \N$ and for $\L$-a.e.~$q\in K$, with derivative $s\mapsto \nabla
w_k(\Phi(q,s))v(\Phi(q,s))$. Passing to a subsequence and relabelling if
necessary, we conclude that
\[
 \lim_{k\to\infty} w_k(\Phi(q,\cdot)) =
   w(\Phi(q,\cdot))
   \quad \mbox{in $W^{1,1}((-\r,\r);\R^2)$}
\]
for $\L$-a.e.~$q \in K$, and that the weak derivative of $w
(\Phi(q,\cdot))$ is $\nabla w (\Phi(q,\cdot)) v(\Phi(q,\cdot))$.

Let $\phi \in C_c^{\infty}((\sigma_q(-\r),\sigma_q(\r));\R^2)$.
Combining the previous observation, \eqref{dif} and \eqref{gioia},
we get
\[
\begin{split}
   \int_{\sigma_q(-\r)}^{\sigma_q(\r)} w(\ga_q(s)) \dot{\phi}(s) \, ds
   & =
   \int_{-\r}^{\r} w(\Phi(q,s)) \dot{\phi}(\sigma_q(s))
              \dot{\sigma}_q(s) \, ds \\
   & =
   -\int_{-\r}^{\r} \nabla w (\Phi(q,s)) v(\Phi(q,s)) \, \phi(\sigma_q (s)) \, ds \\
   & =
   -\int_{\sigma_q(-\r)}^{\sigma_q(\r)} \nabla w (\ga_q(s)) \dot{\ga}_q(s) \, \phi(s) \, ds.
\end{split}
\]
Hence $w \circ \ga_q \in W^{1,1}((\sigma_q(-\r),\sigma_q(\r));\R^2)$,
and its weak derivative is $(\nabla w \circ \ga_q) \dot{\ga_q}$.

\end{proof}

By means of Theorem \ref{regularity}, we can now prove Theorems
\ref{effe} and \ref{acca}.

\medskip

\begin{proof}
[Proof of Theorem \ref{effe}] 
We have $v(z) = 2z-\nabla f^{\perp}(z) \in BV(\O;\R^2) \cap
L^{\infty}(\O;\R^2)$ and $\mathrm{div} \, v = 4$
(cf.~\eqref{coop}). The function $\la:\O \to \R$, $\la = |v|$,
satisfies $0 < c_1 \leq \la \leq c_2$ $\L$-a.e.~in $\O$ by
\eqref{Ok}. By Corollary \ref{reg1}, the vector field $w = v/\la$
belongs to $W^{1,1}(\O;\R^2)$, and $\mathrm{div} \, w^{\perp} = H$
$\L$-a.e.~in $\O$ by \eqref{distr}.

Let $\ga_z:(\sigma_z(-\r),\sigma_z(\r)) \to \O$, be the integral
curve of the vector field $w$ defined in \eqref{gaia} and
satisfying \eqref{gioia}. We claim that $\ga_z$ parameterizes an
arc of circle with curvature $H$. Since $\Phi(z,\cdot)$ is a
reparameterization of $\ga_z$, proving the claim concludes the
proof of Theorem \ref{effe}.

By Remark \ref{compositions with the flow}, the following identities
hold a.e.~in $(\sigma_z(-\r),\sigma_z(\r))$ for $\L$-a.e.~$z \in K$:
\begin{enumerate}
\item
$|w \circ \ga_z| = 1$,
\item
$(w_1 \p_x w_1 + w_2 \p_x w_2) \circ \ga_z
= (w_1 \p_y w_1 + w_2 \p_y w_2) \circ \ga_z = 0$ and
\item
$(\mathrm{div} \, w^{\perp}) \circ \ga_z = H$.
\end{enumerate}
By Theorem \ref{regularity}, we find that
$\ga_z \in W^{2,1}((\sigma_z(-\r),\sigma_z(\r));\R^2)$ with
\begin{equation}\label{ds}
   \ddot{\ga}_z = (\nabla w \circ \ga_z) \dot{\ga}_z.
\end{equation}
Using \eqref{ds}, ii) and iii), we compute
\[
\begin{split}
   \ddot{\ga}_z \cdot \dot{\ga}_z^{\perp}
   =
   (-\mathrm{div} \, w^{\perp}) \circ \ga_z
   =
   - H.
\end{split}
\]
Moreover, by i), we also have $\ddot{\ga}_z \cdot \dot{\ga}_z = 0$
a.e.~in $(\sigma_z(-\r),\sigma_z(\r))$. Then
$\ddot\ga_z= - H \dot\ga_z^\perp$ a.e.~and this implies that $\ga_z$
is of class $C^\infty$. The claim follows.
\end{proof}

\medskip

\begin{proof}[Proof of Theorem \ref{acca}] 
The vector field $v$ in \eqref{uva} is in $BV(\O;\R^2) \cap L^{\infty}(\O;\R^2)$
and satisfies $\mathrm{div} \, v = -4 h_t \in L^{\infty}(\O)$ (cf.~\eqref{migros}).

We denote by $u = (h_y - 2h h_t, 1-2 y h_t)$ the
projection of the vector field \eqref{cosca} onto the $xy$-plane.
The relation between $v$ and $u$ is
\begin{equation}
 \label{rella}
   v = \left( \begin{array}{cc}
     0 & 1 \\
     2y & -2h \\
   \end{array}
   \right) u.
\end{equation}
In \eqref{rella}, we think of $v$ and $u$ as column vectors. The function
$\la:\O \to \R$, $\la = |u|$, satisfies $0 < c_1 \leq \la \leq c_2$
$\L$-a.e.~in $\O$ by \eqref{hurra}. Let $F \in C^{\infty}(\R^2;\R^2)$
be a mapping which satisfies \eqref{F} (with $c_1$ in place of $\d$).
We consider the vector fields $v/\la$ and $w = F \circ u$ in $\O$.
Then $w = u/|u|$ a.e.~in $\O$. We have $w \in W^{1,1}(\O;\R^2)$
by Corollary \ref{reg2} and $v/\la \in W^{1,1}(\O;\R^2)$ by Corollary
\ref{reg3}.

Denote by $\ga_{\zeta}:[\sigma_{\zeta}(-\r),\sigma_{\zeta}(\r)] \to \O$ the
integral curve of $v/\la$ defined in \eqref{gaia}. Then the curve $s \mapsto
w (\ga_{\zeta}(s))$ belongs to
$W^{1,1}((\sigma_{\zeta}(-\r),\sigma_{\zeta}(\r));\R^2)$ for $\L$-a.e.~$\z
\in K$, and its weak derivative is equal to $(\nabla w \circ \ga_\z)\dot{\ga}_{\z}$
(Theorem \ref{regularity}). We claim that
\begin{equation}
  \label{butterfly}
   \nabla w \, \frac{v}{\la} = -H w^{\perp}
   \qquad \textrm{$\L$-a.e.~in $\O$}.
\end{equation}
Indeed, using \eqref{DF}, we compute
\[
\begin{split}
   \nabla  w \, \frac{v}{\la}
   & =
   (\nabla F \circ u) \, \nabla u \, \frac{v}{\la}
   =
   \frac{1}{|u|^4} \big( u^\perp \otimes u^\perp\big) \nabla u\, v \\
   & =
   \frac{1}{|u|^3} \big(u^\perp \cdot (\nabla u \, v) \big) w^\perp.
\end{split}
\]
On the other hand, from \eqref{pox}, we have $\LL w = H$ in $\O$,
where $\LL$ is the differential operator appearing in \eqref{DO} with
$a = 2h$ and $b = -2y$. Using \eqref{DF} and the notation from \eqref{lulu},
we get
\[
\begin{split}
   \LL w
   & =
   \mathrm{tr}\big( (\nabla F \circ u) \nabla u B \big)
   =
   \frac{1}{|u|^3} \mathrm {tr} \big( ( u^\perp \otimes u^\perp ) \nabla u B\big) \\
   & =
   \frac{1}{|u|^3} u^\perp \cdot( u^\perp \nabla u B).
\end{split}
\]
Using \eqref{rella}, by a short computation, we find
\[
   u^\perp \cdot (u^\perp \nabla u B) = -u^\perp \cdot (\nabla u \, v),
\]
and this ends the proof of \eqref{butterfly}. From Remark \ref{compositions
with the flow}, it follows that the identity \eqref{butterfly} holds along
$\L$-a.e.~curve $\ga_\z$. Hence, by Theorem \ref{regularity} applied to
$w$, we have
\begin{equation}
 \label{dirx}
 \frac{d}{d s} (w \circ \ga_\z) = - H (w^{\perp} \circ \ga_\z)
\end{equation}
in the sense of weak derivatives.

The projection of the vector field $\z = (y,t) \mapsto (h(\z),\z)$
onto the $xy$-plane belongs to $W^{1,1}(\O;\R^2)$. Denote by
$\kappa_{\zeta}:[\sigma_{\zeta}(-\r),\sigma_{\zeta}(\r)] \to \R^2$ the
projection of the curve $(h (\ga_{\zeta}),\ga_{\zeta})$ onto the
$xy$-plane. This projection is a reparameterization of the curve in
\eqref{migrol}. We claim that $\kappa_\z$ parameterizes an arc of circle
with curvature $H$ oriented clockwise. By Theorem \ref{regularity},
for $\L$-a.e.~$\z \in K$, we have
$\kappa_\z \in W^{1,1}((\s_{\z}(-\r),\s_{\z}(\r));\R^2)$, and a
short computation shows that the weak derivative of $\kappa_{\z}$ is
$\dot{\kappa}_{\z} = w \circ \ga_{\z}$. From \eqref{dirx}, we deduce
$\kappa_\z \in W^{2,1}((\s_{\z}(-\r),\s_{\z}(\r));\R^2)$ and
$\ddot{\kappa}_{\z} = -H \dot{\kappa}_{\z}^{\perp}$. This implies that
$\kappa_{\z}$ is of class $C^{\infty}$ and the claim follows.
\end{proof}
\section{Characterization of convex isoperimetric sets}
\label{final}
\setcounter{equation}{0}

In this final section, we combine the geometric description of the
characteristic set of a convex set (Theorem \ref{H6}) with the
``foliation by geodesics property'' of the boundary of convex
isoperimetric sets (Corollaries \ref{geodesics1} and
\ref{geodesics2}), in order to prove the Characterization Theorem
\ref{UN}.

\begin{proof}
  [Proof of Theorem \ref{UN}]
According to Theorem \ref{H6}, the characteristic set $\Sigma(C)$
of a convex isoperimetric set $C$ admits a disjoint decomposition
$\Sigma(C) = \Sigma^- \cup \Sigma^+$, where $\Sigma^-$ and
$\Sigma^+$ are two closed, horizontal segments. The curvature
$H$ of $C$ scales like $1/\la$ with respect to the dilations
$\d_\la$, $\la > 0$. Thus, we can assume $H = 2$ without loss
of generality.

\medskip

{\it Claim 1}. For all $p \in \partial C \setminus \Sigma(C)$,
there is a curve of maximal length passing through $p$ and contained
in $\p C$ whose projection onto the $xy$-plane is an arc of circle
with curvature $H$. One endpoint of the curve belongs to $\Sigma^-$,
the other to $\Sigma^+$.

\medskip

Let $f:D \to \R$ be the function in \eqref{dado},  and let
$h:E \to \R$ be the function in \eqref{tratto}. Using an isometry
and a left translation in the $y$ direction if necessary, we can
assume that any point $p \in \p C\setminus\Sigma(C)$ is either of
the form
\begin{itemize}
  \item [1)] $p = (z,f(z))$ for some
          $z \in \mathrm{int}(D)\setminus\Sigma(f)$, or of the
          form
  \item [2)] $p =(h(\zeta),\zeta)$ for some
      $\zeta \in \mathrm{int}(E) \setminus
      \big( \{y=0\}\cup\Sigma(h)\big)$.
\end{itemize}
Then the existence of a piece of geodesic with curvature $H$
passing through $p$ and contained in $\p C$ follows either from
Corollary \ref{geodesics1} or from Corollary \ref{geodesics2}. The
existence of a curve with maximal length whose projection onto the
$xy$-plane is an arc of circle with curvature $H$ follows from a
compactness argument. If one endpoint of the maximal curve does
not belong to $\Sigma(C)$, then a continuation argument involving
the Uniqueness Lemma \ref{uniqueness} provides a proper extension
of the curve, which is impossible. Finally, the endpoints of the
maximal curve cannot both belong to $\Sigma^-$ or $\Sigma^+$. This
follows from our discussion of geodesics in Section \ref{prel}.

\medskip

{\it Claim 2}. The closed horizontal segments $\Sigma^-$ and $\Sigma^+$
are points.

\medskip

Suppose by contradiction that $\Sigma^-$ is a closed, horizontal
segment of positive length. After an isometry and a left
translation, we have $\Sigma^- = \big\{ (0,y,0)\in\H \; \big| \; |y| \leq
y_0 \big\}$ for some $y_0 > 0$. Since the horizontal planes
$\HB_{p_0}$ and $\HB_{-p_0}$ with $p_0=(0,y_0,0)$ are supporting
planes for $C$, we have
\begin{equation}
 \label{hb}
   C\sub \HB_{p_0}^+\cap \HB_{-p_0}^+.
\end{equation}
Claim 1 combined with an approximation-compactness argument shows
that there is a geodesic $\ga:[0,L]\to \p C$ parameterized
proportionally to arc-length, with curvature $H = 2$ and such that
$\ga(0) = 0$. We have $\dot\ga(0) = (\a,\b,0)$ for some $\a, \, \b
\in \R$ with $|(\a,\b)| = 1/2$. Assume $\a = 0$ and $\b = 1/2$ (if
$\a \neq 0$, the proof is easier). Then, by \eqref{seed},
\[
   \ga(s) = \frac 12 \big(1-\cos s,\sin s, s-\sin s\big)\in \p C
   \quad \mbox{for all $s \in [0,L]$}.
\]
However, for any $c_0 > 0$, there is $\d > 0$ with
\[
        s-\sin s < c_0(1-\cos s)\qquad \textrm{for all }s\in
        (0,\d).
\]
This contradicts \eqref{hb}.

\medskip

\it Claim 3. \rm Up to a left translation, the set $C$ is of the
form \eqref{bubble}.

\medskip

After a left translation, we can assume that
$\Sigma^- = \{ (0,-\pi/2) \}$ and that $\Sigma^+ = \{ (z,t) \}$
with $t > -\pi/2$. This forces $z=0$ and $t=\pi/2$, otherwise there
would be at most one geodesic with curvature $H=2$ connecting
$\Sigma^-$ to $\Sigma^+$ and, by Claim 1, this is not possible.

We claim that $0 \in \mathrm{int}(D)$. If, by contradiction, $0
\in \p D$, then $\{ (0,t) \in \H \; | \; |t| \leq \pi/2 \}$ is
contained in $\p C$. But there is no geodesic with curvature
$H = 2$ which starts from $\Sigma^-$ and passes through a point $(0,t)$
with $|t| < \pi/2$. This contradicts Claim 1.

For all  $z \in \p D$, there is a geodesic contained in $\partial C$,
with endpoints $(0,-\pi/2)$, $(0,\pi/2)$, and which passes through $(z,f(z))$.
The projection of this geodesic onto the $xy$-plane is a (full) circle with
radius $1/2$ passing through the points $0$ and $z$. Hence $\partial C$
contains the surface of rotation around the $t$-axis generated by the curve
\eqref{seed}, and the claim follows.
\end{proof}


\end{document}